\documentclass[11pt]{amsart}
\usepackage{amssymb}
\usepackage{amsmath}
\usepackage{amsrefs}
\usepackage[colorlinks=true]{hyperref}
\usepackage{cleveref}
\usepackage{amscd}
\usepackage{amsthm}
\usepackage{mathtools}
\usepackage{tikz-cd}
\usepackage{mathrsfs}
\usepackage{pinlabel}
\usepackage{bm}
\usepackage[margin=2.93cm]{geometry}

\usepackage{xcolor}
\hypersetup{
    colorlinks,
    linkcolor={red!50!black},
    citecolor={green!50!black},
    urlcolor={blue!80!black}
}
\allowdisplaybreaks

\newtheorem{theoremletter}{Theorem}
\renewcommand{\thetheoremletter}{\Alph{theoremletter}}
\newtheorem{theorem}{Theorem}[section]
\newtheorem{lemma}[theorem]{Lemma}
\newtheorem{proposition}[theorem]{Proposition}
\newtheorem{corollary}[theorem]{Corollary}

\newtheorem{conjecture}[theorem]{Conjecture}
\newtheorem{question}[theorem]{Question}

\theoremstyle{definition}
\newtheorem{definition}[theorem]{Definition}
\newtheorem{example}[theorem]{Example}

\theoremstyle{remark}
\newtheorem{remark}[theorem]{Remark}

\numberwithin{equation}{section}

\newcommand{\R}{\mathbb{R}}

\newcommand{\Z}{\mathbb{Z}}

\newcommand{\lt}{\left}
\newcommand{\rt}{\right}
\newcommand{\tms}{\times}

\newcommand{\rmk}{\begin{remark}}
\newcommand{\ermk}{\end{remark}}
\newcommand{\cor}{\begin{corollary}}
\newcommand{\ecor}{\end{corollary}}
\newcommand{\eq}{\begin{equation}}
\newcommand{\eeq}{\end{equation}}
\newcommand{\eqs}{\begin{equation*}}
\newcommand{\eeqs}{\end{equation*}}
\newcommand{\prop}{\begin{proposition}}
\newcommand{\eprop}{\end{proposition}}
\newcommand{\thm}{\begin{theorem}}
\newcommand{\ethm}{\end{theorem}}
\newcommand{\conj}{\begin{conjecture}}
\newcommand{\econj}{\end{conjecture}}
\newcommand{\lem}{\begin{lemma}}
\newcommand{\elem}{\end{lemma}}
\newcommand{\defi}{\begin{definition}}
\newcommand{\edefi}{\end{definition}}
\newcommand{\ex}{\begin{example}}
\newcommand{\eex}{\end{example}}
\newcommand{\alis}{\begin{align*}}
\newcommand{\ealis}{\end{align*}}
\newcommand{\pf}{\begin{proof}}
\newcommand{\epf}{\end{proof}}
\newcommand{\ali}{\begin{align}}
\newcommand{\eali}{\end{align}}
\newcommand{\qus}{\begin{question}}
\newcommand{\equs}{\end{question}}

\newcommand{\mc}{\mathcal}
\renewcommand{\bf}{\textbf}

\newcommand{\C}{\mathbb{C}}

\newcommand{\sub}{\subset}

\newcommand{\ov}{\overline}

\newcommand{\bb}{\mathbb}

\newcommand{\op}{\operatorname}

\renewcommand{\a}{\alpha}
\renewcommand{\b}{\beta}
\renewcommand{\d}{\partial}
\newcommand{\e}{\epsilon}

\newcommand{\g}{\gamma}

\renewcommand{\t}{\theta}

\renewcommand{\o}{\omega}

\newcommand{\fk}{\frak}

\newcommand{\G}{\Gamma}

\renewcommand{\S}{\Sigma}

\renewcommand{\ov}{\overline}

\newcommand{\HP}{\mathit{HP}}

\begin{document}
\title[On the sheaf-theoretic  $\operatorname{SL}(2,\mathbb{C})$ Casson-Lin invariant]{On the sheaf-theoretic  $\mathbf{SL}(2,\mathbb{C})$ Casson-Lin invariant}
\author[Laurent C\^ot\'e]{Laurent C\^ot\'e}
\thanks {LC was supported by a Stanford University Benchmark Graduate Fellowship.}
\address{Department of Mathematics, Stanford University, 450 Serra Mall, Stanford, CA 94305}
\email{lcote@stanford.edu}

\author[Ikshu Neithalath]{Ikshu Neithalath}
\address {Department of Mathematics, UCLA, 520 Portola Plaza, Los Angeles, CA 90095}
\email {ikshu@math.ucla.edu}

\begin{abstract} We prove that the ($\tau$-weighted, sheaf-theoretic) $\op{SL}(2,\C)$ Casson-Lin invariant introduced by Manolescu and the first author in \cite{cote-man} is generically independent of the parameter $\tau$ and additive under connected sums of knots in integral homology $3$-spheres. This addresses two questions asked in \cite{cote-man}.  Our arguments involve a mix of topology, microlocal analysis and algebraic geometry, and rely crucially on the fact that the $\op{SL}(2,\C)$ Casson-Lin invariant admits an alternative interpretation via the theory of Behrend functions.
\end{abstract}
\maketitle

\section{Introduction}  

\subsection{Statement of results}

Given a knot $K$ in a closed, orientable $3$-manifold $Y$ and a real parameter $\tau \in (-2,2)$, Manolescu and the first author defined in \cite{cote-man} a sequence of abelian groups $\HP^*_{\tau}(K)$ which are topological invariants. These groups are constructed using tools from derived algebraic geometry, but they can morally be interpreted as the Morse homology of the $\op{SL}(2,\C)$ Chern-Simons action functional on $Y-K$, restricted to the space of connections with trace $\tau \in (-2,2)$ along the knot meridian. Their precise definition is reviewed in \Cref{subsection:construction}.

The Euler characteristic $ \chi_{\tau}(K):= \sum_{n \in \Z} (-1)^n \op{rk}_{\Z} \HP^n_{\tau}(K)$ is of independent interest since it can be viewed as an $\op{SL}(2,\C)$ analog of the Casson-Lin invariant. The Casson-Lin invariant, which is defined using gauge theory, counts $\op{SU}(2)$ connections with trace zero along the knot meridian and has been well-studied in the literature; cf. \cite{lin-casson, herald}. 

This main goal of this paper is to establish new properties of the $\tau$-weighted $\op{SL}(2,\C)$ Casson-Lin invariant $\chi_{\tau}(-)$, which arguably make it a better invariant from the perspective of topology than its categorification $\HP^*_{\tau}(-)$. In doing so, we also partly answer some questions which were stated in \cite{cote-man}.  

Our first result states that the $\op{SL}(2,\C)$ Casson-Lin invariant is generically independent of the parameter $\tau$.

\begin{theoremletter} \label{theorem:cassonconstant}
	Let $K$ be an oriented knot in a closed, oriented $3$-manifold $Y$. Then the (sheaf-theoretic, $\tau$-weighted) $\op{SL}(2,\C)$ Casson-Lin invariant $\chi_{\tau}(K)$ is constant as a function of $\tau$ on a Zariski open subset of the complex plane.
\end{theoremletter}

\Cref{theorem:cassonconstant} answers a weaker form of Question 1.5 in \cite{cote-man}, which asked whether $\HP^*_{\tau}(-)$ is generically independent of $\tau$. The statement of this theorem merits some clarification due to the fact that $\HP^*_{\tau}(-)$ was only defined in \cite{cote-man} for $\tau \in (-2,2)$. In fact, we will show in \Cref{section:complexplane} that the construction of $\HP^*_{\tau}(-)$ can be generalized to all $\tau \in \C- \{ \pm 2\}$. Moreover, there is an alternative definition of $\chi_{\tau}(-)$ which makes sense for all $\tau \in \C$ and agrees with the Euler characteristic of $\HP^*_{\tau}(-)$ when these groups are defined. Our proof of \Cref{theorem:cassonconstant} actually works with this alternative definition of $\chi_{\tau}(-)$. 

As a result of \Cref{theorem:cassonconstant}, we can introduce the following definition.

\defi \label{definition:generalcasson}
	Let $ \chi_{CL}(K) \in \Z$ be defined as the generic value of $\chi_{\tau}(K)$ for $\tau \in \C$. We say that $\chi_{CL}(-)$ is the \emph{(sheaf-theoretic) $\op{SL}(2,\C)$ Casson-Lin invariant}.
\edefi

Although $\chi_{CL}(-)$ contains less information than $\HP^*_{\tau}(-)$, our next theorem shows that it has the advantage of being additive under connected sums of knots.

\begin{theoremletter} \label{theorem:connectedsum}
	For $i=1,2$, let $Y_i$ be a closed, orientable integral homology $3$-sphere and let $K_i \sub  Y_i$ be a knot. Letting $K_1 \# K_2$ denote the connected sum of $K_1$ and $K_2$ (see \cite[Sec.\ 7.1]{cote-man}), we have $$\chi_{CL}(K_1 \# K_2) = \chi_{CL}(K_1) + \chi_{CL}(K_2).$$ 
\end{theoremletter}

\Cref{theorem:connectedsum} affirmatively answers Question 1.6 of \cite{cote-man} for generic $\tau \in \C$. This question asks whether $\chi_{\tau}(-)$ is additive for knots in $S^3$ and all $\tau \in (-2,2)$, so there are always finitely many cases for which it remains open. However, one could reasonably argue that $\HP^*_{\tau}(-)$ is not a meaningful invariant for certain non-generic choices of $\tau$, and that \Cref{theorem:connectedsum} therefore addresses the most interesting part of Question 1.6. This is because $\HP^*_{\tau}(-)$ only counts irreducible representations, and families of irreducibles can sometimes converge to a reducible representation at certain exceptional points. If $\rho_{red}: \pi_1(Y-K) \to \op{SL}(2,\C)$ is such a representation and has trace $\tau_0 \in \C - \{ \pm 2\}$ along the meridian of $K$, then $\HP^*_{\tau_0}(K)$ does not see $\rho_{red}$ and therefore gives the ``wrong" count. This situation could hopefully be corrected by defining an invariant which also takes into account reducibles.  

We remark that a weaker version of \Cref{theorem:connectedsum} was proved by Manolescu and the first author in \cite[Thm.\ 7.17]{cote-man}. They showed, for $K_i \sub Y_i$ a knot in an integral homology $3$-sphere, that $\chi_{\tau}(\#_{i=1}^n K_i)= \sum_{i=1}^n \chi_{\tau}(K_i)$ for generic $\tau \in \C$ under the assumption that the character schemes $\mathscr{X}^{\tau}_{irr}(K_i)$ are smooth. In principle, \Cref{theorem:connectedsum} is a much stronger result since character schemes of knot complements can be singular in general (in fact, singularities of $3$-manifold groups can in some sense be arbitrarily bad; see \cite{kapovich-millson}).  On the other hand, from a purely computational perspective, \Cref{theorem:connectedsum} may turn out not to be particularly useful: for a knot $K \sub Y$, one needs to understand $\mathscr{X}^{\tau}_{irr}(K)$ very well in order to compute $\chi_{\tau}(K)$. The only examples that the authors have been able to handle turn out to be smooth. 

\subsection{$\op{SL}(2,\C)$ Floer homology in families}
It is natural to ask about the behavior of the invariants $\HP^*_{\tau}(-)$ in families. Optimistically, there should exist for any knot $K \subset Y$ a constructible sheaf $\mc{F}(K) \in D^b(\C- \{\pm 2\} )$ whose stalk at a point $\tau \in \C$ is precisely $\HP^*_{\tau}(K)$. In particular, the groups $\HP^*_{\tau}(K)$ should form a local system on some Zariski open subset of $\C$.  The monodromy of this local system may then carry additional information which is not seen by the stalks.  One possibly naive justification for this picture is that the constructions of \cite{cote-man} mostly take place in the algebraic category, and so one expects that any exceptional behavior should only occur on a Zariski closed subset of the complex plane. A more systematic potential approach is briefly outlined in \Cref{subsection:construction}.  

Observe that if one could show that $\HP^*_{\tau}(K)$ forms a local system on a Zariski open subset of $\C$, then this would imply \Cref{theorem:cassonconstant}. Thus, \Cref{theorem:cassonconstant} could be seen as providing evidence for the validity of the above picture. As mentioned above, we prove in \Cref{section:complexplane} that the construction of the invariants $\HP^*_{\tau}(K)$, which is carried out in \cite{cote-man} for $\tau \in (-2,2)$, actually works for all $\tau \in \C - \{ \pm 2\}$. By combining this with the arguments from Section 5 of \cite{cote-man}, one can easily construct a local system on a Zariski open subset of the complex plane whose stalks compute $\HP^*_{\tau}(K)$ for a wide class of knots whose character variety is one-dimensional; see \Cref{subsection:families}. We explicitly describe this local system for the figure-eight knot in \Cref{subsection:families} and find that it exhibits non-trivial monodromy. 

\ex Let $K=4_1$ denote the figure-eight knot. There is a local system (viewed as an object of the derived category) $\mc{F}(K) \in D^b(\C - \{ \pm 1, \pm \sqrt{5} \})$ whose stalks compute $\HP^*_{\tau}(K)$ in the sense that $H^*(\mc{F}(K)_{\tau}) =\HP^*_{\tau}(K)$. This local system has stalks (non-canonically) isomorphic to $\Z^2$ in degree zero. Fixing an isomorphism over some basepoint, the monodromy map around each puncture has eigenvalues $ i$ and $-i$.  \eex

\subsection{Summary of the proofs}

This paper builds on work of the first author and Manolescu in \cite{cote-man}. Let us therefore try to situate our work with respect to theirs. 

The key new player in this paper is the so-called \emph{Behrend function}, introduced by Behrend \cite{behrend}. The Behrend function is a constructible function which can be associated to any scheme (or complex-analytic space) over $\C$. Roughly speaking, it keeps track of singular or non-reduced behavior on the scheme; in particular, it is identically equal to $(-1)^n$ on the locus of smooth points of an $n$-dimensional scheme. The precise definition is explained in \Cref{section:behrend-function}.

Given a knot $K \sub Y$, the $\op{SL}(2,\C)$ Casson-Lin invariant $\chi_{\tau}(K)$ can be defined in two ways. The first definition is just the one above, namely the alternating sum of the ranks of the groups $\HP^*_{\tau}(K)$ constructed in \cite{cote-man}. The second one defines $\chi_{\tau}(K)$ as the Euler characteristic of the character scheme $\mathscr{X}^{\tau}_{irr}(K)$ \emph{weighted} by the Behrend function (see \Cref{definition:weightedeuler}). The fact that these definitions agree is essentially built into Joyce's theory of critical loci, which underlies the construction of the groups $\HP^*_{\tau}(K)$ in \cite{cote-man}. 

The perspective taken in this paper is to work almost entirely with the second definition, in terms of Behrend functions. This is in contrast to \cite{cote-man}, which only considers the first definition.

The usefulness of this perspective is illustrated by our proof of \Cref{theorem:cassonconstant}, Indeed, \Cref{theorem:cassonconstant} is deduced as a straightforward corollary of the following result, which is purely a statement about complex algebraic geometry and may be of independent interest.

\renewcommand\thetheoremletter{A'}

\begin{theoremletter} \label{theorem:a1} Fix a complex affine variety $X$ and a morphism $X \to \bb{A}^1$. Then the function $\tau \mapsto \chi_B(X_{\tau})$ is Zariski-locally constant (here $\chi_B(-)$ denotes the Euler characteristic weighted by the Behrend function).
\end{theoremletter}

The statement of \Cref{theorem:a1} appears plausible in light of the general philosophy that ``bad behavior should only occur in codimension $\geq 1$". However, the proof is in our opinion nontrivial, and crucially exploits deep work of Verdier \cite{verdier} on stratifications of complex algebraic varieties.

Our proof of \Cref{theorem:connectedsum} also relies on $\chi_{\tau}(K)$ being defined via the Behrend function. As mentioned above, \Cref{theorem:connectedsum} generalizes Theorem 7.17 in \cite{cote-man}. The proof of \Cref{theorem:connectedsum} can in fact be carried out along similar lines as the original argument in \cite{cote-man}, treating the Euler characteristic weighted by the Behrend function as a replacement for the topological Euler characteristic which was considered in the proof of Theorem 7.17 in \cite{cote-man}. However, there were some technical challenges to this generalization: the proof of \Cref{corollary:typeiizero} relies on nontrivial results about triangulations of algebraic subsets, and we also needed to clarify some ambiguities in the literature concerning the definition of the Behrend function, which is done in the appendix.  

In \Cref{section:complexplane}, we prove that the invariant $\HP^*_{\tau}(K)$ constructed for $\tau \in (-2,2)$ in \cite{cote-man} can be defined for all $\tau \in \C- \{ \pm 2\}$. This actually boils down to proving that a certain character variety $X^{\tau}_{irr}(\S)$ is connected and simply-connected. For $\tau \in (-2,2)$, this character variety is homeomorphic to a certain moduli space of Higgs bundles. This homeomorphism was exploited in \cite{cote-man} to analyze their topology. 

For $\tau \in \C-\{\pm 2 \}$, one needs to consider so-called $K(D)$-pairs, which are a mild generalization of Higgs bundles. Using again a result of Verdier on stratifications of complex varieties \cite{verdier} and a ``variation of weights" argument originally due to Thaddeus \cite{thaddeus} which also played a role in \cite{cote-man}, we deduce that $X^{\tau}_{irr}(\S)$ is connected and simply connected for all $\tau \in \C-\{\pm 2\}$ as a consequence of the fact, proved in \cite{cote-man}, that this holds for $\tau \in (-2,2)$. 

This paper is intended to be accessible to topologists who are familiar with the language of schemes. We have therefore tried to give detailed arguments and references to the algebraic-geometry literature, some of which would presumably be unnecessary in a paper aimed only at algebraic geometers.

\subsection*{Acknowledgements} 
We thank Brian Conrad, Tony Feng, Dominic Joyce, Aaron Landesman, Ciprian Manolescu, Lisa Sauermann, Vivek Shende and Burt Totaro for many helpful conversations and suggestions. 

\section{Context and background material}

\subsection{Review of the construction of $\HP^*_{\tau}(-)$ in \cite{cote-man}} \label{subsection:construction}

Fix the data of an oriented knot $K$ in a closed, orientable $3$-manifold $Y$ and a real parameter $\tau \in (-2,2)$.  Let $E_K$ denote the knot exterior and fix a suitable Heegaard splitting $(\S, U_0, U_1)$ of $E_K$, where $\S$ is a Riemann surface with two disks removed and the $U_i$ are handlebodies. Now consider the moduli space $X^{\tau}_{irr}(\S)$ of irreducible $\op{SL}(2,\C)$ representations with trace $\tau$ along the boundary circles. This moduli space turns out to admit a natural holomorphic symplectic structure and the natural inclusions $\iota_j: L_j:= X^{\tau}_{irr}(U_j) \hookrightarrow X^{\tau}_{irr}(\S)$ are Lagrangian embeddings.

Work of Bussi \cite{bussi} based on Joyce's theory of d-critical loci then gives a perverse sheaf $P^{\bullet}_{L_0, L_1}$ on (the complex-analytification of) $X^{\tau}_{irr}(E_K)$. This perverse sheaf is built by locally writing $L_1$ as the critical locus of a holomorphic function $f$ on a small open $\mc{U} \sub L_0$ and applying the vanishing cycles functor to the constant sheaf $\Z|_{\mc{U}}$ to get a perverse sheaf on $\op{crit}(f)= \mc{U} \cap L_0 \cap L_1$. If $L_0, L_1$ are equipped with spin structures, then after tensoring these locally defined perverse sheaves with an appropriate line bundle which keeps track of spin structures on $L_0, L_1$, one shows that the local data can be patched together to give a perverse sheaf on $L_0 \cap L_1$. It can be shown that this perverse sheaf is independent of the Heegaard splitting (see \cite[Prop.\ 3.9]{cote-man}), and $\HP^*_{\tau}(K)$ is defined to be its hypercohomology.

One can attempt to carry out this construction in families, and this could be relevant to the proposal discussed in the introduction that the groups $\HP^*_{\tau}(K)$ should form a constructible sheaf on $\C- \{ \pm 2\}$.  To explain this, let us denote by $X_{irr}(\S)$ the moduli space of irreducible, flat $\op{SL}(2,\C)$ connections on $\S$ having the same holonomy up to conjugacy along the two boundary circles. There is a map $\pi: X_{irr}(\S) \to \C$ taking a representation to the trace of its holonomy along the boundary circles, and the fibers are precisely $X^{\tau}_{irr}(\S)$. Let $\{L_0^{\tau}\}, \{L_1^{\tau}\}$ be the family of Lagragians associated to a Heegaard splitting as above. One way to produce a constructible sheaf on $\C- \{ \pm 2\}$ whose stalks compute $\HP^*_{\tau}(K)$ is to exhibit a constructible sheaf $\mc{F}$ on $\bigcup_{\tau} (L_0^{\tau} \cap L_1^{\tau})$ such that $i_{\tau}^! \mc{F}= P^{\bullet}_{\tau}(K)$, where $i_{\tau}: X_{irr}^{\tau}(\S) \to X_{irr}(\S)$ is the natural inclusion. Indeed, the proper base change theorem \cite[Thm.\ 3.2.13]{dimca} then implies that $H^*((\pi_! \mc{F})_{\tau}) = \bb{H}^*(i_{\tau}^! \mc{F})$. 

To construct such a sheaf, one approach could be to locally write $\bigcup_{\tau} (L_0^{\tau} \cap L_1^{\tau})$ as the critical locus of a holomorphic family of holomorphic functions $f_{\tau} (z_1, \dots, z_n)$ on some open $\mc{U} \sub X_{irr}(\S)$.  One can apply the vanishing cycles functor to this locally defined family which should produce a constructible sheaf $\mc{F}_U$ on $\bigcup_{\tau}(L_0^{\tau} \cap L_1^{\tau}) \cap \mc{U}$. One can then ask whether, after tensoring with suitable spin data, these constructible sheaves can be patched together.

\subsection{Notation and conventions}  As a general rule, we always use the same notation and follow the same conventions as in \cite{cote-man}. 

\begin{itemize}
\item All schemes are assumed to be separated and of finite-type over $\C$. 
Since we are exclusively dealing with subschemes of affine varieties, these hypotheses are automatically satisfied.  
\item  All subschemes are assumed to be locally closed. 
\item As in \cite{cote-man}, a variety is a (not necessarily irreducible) reduced scheme. A subvariety of a scheme is a subscheme which is a variety. 
\item Following \cite[Sec.\ 1.3]{fulton}, an \emph{algebraic cycle} on a scheme $X$ is a finite formal sum of irreducible subvarieties with integer coefficients. These form a group under addition which is denoted $Z_*(X)$. 
\item When we refer to a point of a $\C$-scheme, we mean a closed point unless otherwise indicated. To lighten the notation, we will not distinguish between $\C$-schemes and their associated set of closed points $X(\C)$ in situations where the intended meaning seems clear.  Thus, if $X$ is a subscheme of $Y$, we sometimes write $X \sub Y$ as shorthand for $X(\C) \sub Y(\C)$. 
\item If $X_1,\dots,X_n$ are subschemes of $X$, then $\sqcup X_i \sub X$ is the \emph{set-theoretic} union of the points of $X_i$. In particular, $\sqcup X_i$ should be viewed as a topological subspace of $X$ endowed with the subspace topology -- not with the disjoint union topology. 
\item Given a scheme $X$, a \emph{partition} is a collection of pairwise disjoint subschemes $\{X_i\}$ such that $X= \sqcup X_i$.
\end{itemize}

We now give a brief overview of some constructions and objects described in more detail in \cite{cote-man} and which we will also be using. 

We will be considering Heegaard splittings $U_0 \cup_{\S} U_1$ of a knot complement $Y-K$ or knot exterior $E_K$. In case $U_0 \cup_{\S} U_1= Y-K$, we always assume as in \cite{cote-man} that $\S$ has genus at least six, that $K$ intersects $\ov{\S}$ in two points, and that the arcs $K \cap \ov{U_i}$ are isotopic rel endpoints to arcs contained in $\ov{\S} = \d \ov{U_i}$. There are analogous conditions for Heegaard splittings of knot exteriors; cf.\ \cite[Def.\ 3.2.]{cote-man}. Concretely, such Heegaard splittings can be constructed by choosing a Morse function on $Y$ with a single minimum and maximum on $K$, and such that $K$ is preserved by the gradient flow for an auxiliary metric.

Given a finitely-presented group $\G = \langle g_1,\dots,g_k \mid r_1,\dots, r_l \rangle$, we let $\mathscr{R}(\G)$ be its $\op{SL}(2,\C)$ \emph{representation scheme}. This is an affine scheme over $\C$ whose closed points parametrize $\op{SL}(2,\C)$ representations of $\G$. Let $R(\G) \sub \mathscr{R}(\G)$ be the unique reduced closed subscheme with the same topology which we call the \emph{representation variety}. The group scheme $\op{SL}_2$ acts by conjugation on $\mathscr{R}(\G)$ and the GIT quotient is called the \emph{character scheme} and denoted by $\mathscr{X}(\G)$. We let $X(\G) \sub \mathscr{X}(\G)$ be the associated reduced subscheme and refer to it as the \emph{character variety}.

We will always be considering representation schemes/varieties arising from Heegaard splittings of knot complements or exteriors. To lighten the notation, we write $\mathscr{R}(\S)= \mathscr{R}(\pi_1(\S)),$ $\mathscr{R}(U_i)= \mathscr{R}(\pi_1(U_i)), \mathscr{R}(K)= \mathscr{R}(Y-K)= \mathscr{R}(\pi_1(Y-K))$, where $Y= U_0 \cup_{\S} U_1$ and the fundamental group is always assumed to be defined with respect to a basepoint on $\S$. We use analogous shorthand notation for representation varieties and for character varieties and schemes, and for Heegaard splittings of knot exteriors. 

Finally, there are \emph{relative} versions of the above notions. More precisely, given a finitely presented group $\G$ and some conjugacy classes $\fk{c}_1,\dots, \fk{c}_j \sub \G$, one can consider relative representation schemes which parametrize representations with fixed trace $\tau$ on the $\fk{c}_i$.  These behave essentially like ordinary representation schemes, and a detailed account is provided in \cite[Sec.\ 2.1]{cote-man}.  For $Y-K= U_0 \cup_{\S} U_1 $, let $\mathscr{R}^{\tau}(\S)$ be the relative representation scheme parametrizing representations with fixed trace along the two boundary punctures.  Let $\mathscr{R}^{\tau}(U_i)$ be the scheme of representations having fixed trace along the knot meridian. Let $\mathscr{R}^{\tau}(K)= \mathscr{R}^{\tau}(E_K) = \mathscr{R}^{\tau}(Y-K)$ be the scheme of representations of $ \pi_1(Y-K)=\pi_1(E_K)$ having fixed trace along the knot meridian. We use analogous notation to denote relative representation varieties, and to denote relative character schemes and varieties. 

All of the schemes and varieties described above have an open locus consisting of irreducible representations. We denote them by $\mathscr{R}_{irr}(\G) \sub \mathscr{R}(\G)$, and similarly for the other cases.  We remark that the character varieties $X(\S), X(U_i), X(K)$ (resp. $X^{\tau}(\S), X^{\tau}(U_i), X^{\tau}(K)$) can also be viewed as moduli spaces of flat $\op{SL}(2,\C)$ connections (resp. flat connections with holonomy having trace $\tau$ along the knot meridian). The arguments of this paper do not rely on this interpretation, but we have used it informally in the introduction.

\subsection{Stratifications and constructible objects} \label{subsection:constructible-theory}

Given two vector subspaces $F, G \sub \R^n$, we let $$\delta(F, G):= \sup_{\substack{x \in F\\ \|x\|=1}} \op{dist}(x, G),$$ where $\| \cdot \|$ is the standard Euclidean metric and the distance is also measured using this metric.

Let $M$ and $M'$ be smooth, locally closed submanifolds of $\R^n$ such that $M \cap M'= \emptyset$ and $y \in \ov{M} \cap M'$.  Following \cite[Sec.\ 1]{verdier}, we say that the pair $(M, M')$ satisfies property w) at $y$ if there is a neighborhood $U$ of $y$ in $\R^n$ and a positive constant $C$ such that for all $q' \in U \cap M'$ and $x \in U \cap M$ we have $\delta(T_{M', y'}, T_{M, x}) \leq C \| x -y' \|$. We say that the pair $(M, M')$ satisfies property w) if it satisfies this property at all points $y \in \ov{M} \cap M'$. 

Let $M, M'$ be locally closed submanifolds of a complex algebraic variety $V$ such that $M \cap M'= \emptyset$ and $y \in \ov{M} \cap M'$. We say that the pair $(M, M')$ satisfies condition w) at $y$ if there is a local real analytic embedding $\phi: U \cap V \to \R^n$ so that $(\phi(M), \phi(M'))$ satisfies property w) at $\phi(y)$. We say that $(M, M')$ satisfies property w) if it satisfies this condition for all $y \in \ov{M} \cap M$. 

We now introduce the notion of a w-stratification. 

\defi[see (2.1) in \cite{verdier}]  \label{definition:wstrat} A \emph{w-stratification} of a variety $X$ (i.e. a reduced scheme) over $\C$ is a partition $X= \sqcup_{i=1}^n X_i$, where the $X_i \sub X$ are smooth, connected subschemes, which satisfies the following axioms:
\begin{itemize}
\item[(i)] $X_i \cap X_j= \emptyset$ if $i \neq j$. 
\item[(ii)] If $\ov{X}_i \cap X_j \neq \emptyset$, then $X_j \sub \ov{X}_i$. (One gets the same notion using the analytic or Zariski topology.)
\item[(iii)] If $X_i \sub \ov{X}_j$ and $i \neq j$, then the pair $(X_j, X_i)$ satisfies the condition w). 
\end{itemize}
A w-stratification of a $\C$-scheme just means a w-stratification of the associated variety. The notion of a w-stratification is introduced by Verdier in \cite[(2.1)]{verdier}.  Unless otherwise specified, we only consider w-stratifications in this paper. We will therefore usually omit the prefix and refer to w-stratifications simply as \emph{stratifications}. \edefi

\rmk It is shown in \cite{verdier} that w-stratifications are Whitney stratifications (i.e. they satisfy Whitney's so-called (b) condition). The converse is in general not true. We have chosen to work with w-stratifications in this paper simply for consistency with \cite{verdier} since we quote results of this paper throughout. However, we don't use any properties of w-stratifications which aren't also satisfied by Whitney stratifications, so we could just as easily have worked with ordinary Whitney stratifications. \ermk 

\defi \label{definition:constructible-function} Given a scheme $X$, a subset $C \sub X(\C)$ is said to be \emph{constructible} if it is a finite union of subschemes, i.e. $C= \sqcup_{i=1}^nX_i(\C)$, where $X_i(\C) \sub X(\C)$ is a subscheme. A function $f: X(\C) \to \Z$ is said to be \emph{constructible} if $f(X)$ is finite and if, for every $n \in \Z$, $f^{-1}(n) \sub X(\C)$ is constructible; see \cite{joyce-jlms}. \edefi

It follows from \cite[(2.2)]{verdier} that we can always refine our partition to be a stratification; in particular, we can assume that the $X_i$ are smooth. 

\defi[see \cite{behrend} or Sec.\ 3.3 of \cite{jiang-thomas}] \label{definition:weightedeuler} Let $f: X \to \Z$ be a constructible function. We define the Euler characteristic of $X$ \emph{weighted} by the constructible function $f$ as \eq \chi(X, f):= \sum_{n \in \Z} n \chi(f^{-1}(n)), \eeq where $\chi(-)$ is the topological Euler characteristic.\edefi

We warn the reader that there is considerable ambiguity in the literature concerning the definition of the Euler characteristic weighted by a constructible function. First of all, some authors define $\chi(X,f)$ as in \Cref{definition:weightedeuler} but replacing the topological Euler characteristic with the Euler characteristic with compact support. Other authors adopt the following definition. Given a constructible set $C \sub X(\C)$ and a partition $C= \sqcup_{i=1}^m X_i(\C)$ where the $X_i \sub X$ are subschemes, one defines $\chi_{an}(C):= \sum_{i=1}^n \chi_{c}(X_i(\C))$, where $\chi_{c}(X_i(\C)):= \sum_{k \in \Z} (-1)^k H^k_c(X_i(\C))$ is the Euler characteristic with compact support of $X_i(\C)$; see \cite[Def.\ 3.7]{joyce-jlms}. One can then show (see \cite[(2)]{joyce-jlms}) that $\chi_{an}(C)$ is independent of the chosen partition. One then sets $\chi (X,f):= \sum_{n \in \Z} n \chi_{an}(f^{-1}(c))$. 

It is not at all obvious that these definitions all agree in the present context. We have therefore provided a proof of their equivalence in the appendix. We wish to emphasize that the argument in the appendix uses the fact that we are dealing with varieties over $\C$ (the analogous equivalences would be false over $\R$).

\subsection{The Behrend function}  \label{section:behrend-function}

Given a $\C$-scheme $X$, there exists a constructible function $\nu_X: X(\C) \to \Z$ introduced by Behrend in \cite{behrend} which is usually called the \emph{Behrend function}; see also \cite[Sec.\ 4.1]{joyce-song}. Let us state the definition for affine subschemes of $\bb{A}^n$ as this is sufficient for our purposes. 

Let $X \sub \bb{A}^n$ be an affine scheme over $\C$ defined by the ideal $I \sub \C[x_1, \dots, x_n]$.  We can consider the $\C$-algebra $R= \bigoplus_{m \geq 0} I^m/I^{m+1}$ (where $I^0:= \C[x_1,\dots,x_n]$) and let $$ C_{X/\bb{A}^n}=  \op{Spec} R.$$ The $\C$-algebra inclusion $\C[x_1,\dots, x_n]/ I \hookrightarrow R$ induces a projection map $\pi: C_{X/ \bb{A}^n} \to X$. We say that $C_{X/\bb{A}^n}$ is the \emph{normal cone} of $X \sub \bb{A}^n$; see \cite[B.6]{fulton}. This is a generalization of the normal bundle (and coincides with it for smooth schemes). 

We define \eq \label{equation:canonicalcycle} \fk{c}_{X / \bb{A}^n}:= \sum_{C'} (-1)^{\op{dim} \pi(C')} \op{mult}(C') \pi(C') \in Z_*(X),\eeq where the sum is over all irreducible components $C' \sub C_{X/ \bb{A}^n}$; see \cite[Sec.\ 4.1]{joyce-song}. Here $\pi(C')$ denotes the underlying reduced closed subscheme which is the image of $C'$ under $\pi$. The multiplicity $\op{mult}(C')$ is the length of $C_{X/\bb{A}^n}$ at the generic point of $C'$. This is often referred to as the geometric multiplicity, for instance in \cite[Sec.\ 1.5]{fulton}.  

It turns out that the cycle $\fk{c}_{X/\bb{A}^n}$ depends only on $X$, i.e. it is independent of the embedding of $X$ into $\bb{A}^n$. Letting $\op{CF}(X)$ denote the group of constructible functions on $X$, there is a well-known group morphism $\op{Eu}: Z_*(X) \to \op{CF}(X)$ called the \emph{(local) Euler obstruction} which was originally introduced by MacPherson in \cite{macpherson}. We now define the Behrend function $\nu_X: X  \to \Z$ by letting $\nu_X := \op{Eu}(\fk{c}_{X/\bb{A}^n})$.

If $\nu_X: X \to \Z$ is the Behrend function, then we write $$\chi_B(X):= \chi(X, \nu_X).$$ We refer to this quantity the Euler characteristic of $X$ weighted by the Behrend function.

The Behrend function plays an essential role in this work.  This is mainly due to the fact that it can be computed in two ways on the schemes which we will be considering. On the one hand, the Behrend function can be defined for any finite-type $\C$-scheme in terms of normal cones and Euler obstructions, as explained above for affine schemes.  On the other hand, if we consider a $\C$-scheme $X$ whose complex-analytification is locally the critical locus of a holomorphic function $f: V \sub \C^n \to \C$, then for $x \in V \cap \op{crit}(f) \hookrightarrow X$, the Behrend function can be computed in terms of the Euler characteristic of the sheaf of vanishing cycles of $f$.

More precisely, we have the following formula, due to Parusi\'{n}ski-Pragacz \cite[Thm.\ 4.7]{joyce-song}: \eq \nu_X(x)= (-1)^{\op{dim} V}(1-\chi(\op{MF}_{f}(x))),\eeq where $\op{MF}_f(-)$ is the Milnor fiber of $f$ and $x \in V \cap \op{crit}(f) \hookrightarrow X$. 

Now, let us consider the perverse sheaf $P^{\bullet}_{\tau}(K)$ introduced in \cite{cote-man}, for a knot $K \sub Y$. Recall that $P^{\bullet}_{\tau}(K)$ is a perverse sheaf on (the complex analytification of) the relative chararacter scheme $\mathscr{X}^{\tau}_{irr}(K)$ of $K$, for some $\tau \in (-2,2)$. Recall from the discussion in \Cref{subsection:construction} (see also \cite[p.\ 17]{abou-man}) that a choice of Heegaard splitting $(\S, U_0, U_1)$ allows one to write locally $\mathscr{X}_{irr}^{\tau}(K)$ as the critical locus of a holomorphic function $f: V \to \C$ for some neighborhood $V \sub X^{\tau}(U_0)$.  In this case, we have that $P^{\bullet}_{\tau}(K)|_{V}$ is isomorphic to the perverse sheaf of vanishing cycles of $f$. Hence, one can compute as in \cite[p.\ 20]{abou-man} that, for $x \in V \cap \op{crit}(f) \hookrightarrow  \mathscr{X}_{irr}^{\tau}(K)$, we have
\begin{align*}  \sum_{i \in \Z} (-1)^i \op{rk} \mc{H}^i(P_{\tau}^{\bullet}(K))_x &= \sum_{i \in \Z} (-1)^i \op{rk} \mc{H}^i ( \phi[-1](\Z[\op{dim} V ]) )_x\\
&= (-1)^{\op{dim} V -1}(\chi(\op{MF}_f(x))-1) \\
&= (-1)^{\op{dim} V }(1-\chi(\op{MF}_{f}(x))).
\end{align*}

According to \cite[Thm.\ 4.1.22]{dimca}, if $\mc{F}^{\bullet}$ is a complex of sheaves constructible with respect to a Whitney stratification $\mc{S}$, then 
\eq \sum_{i \in \Z} (-1)^i \op{rk} \bb{H}^i(X, \mc{F}^{\bullet})= \chi(X, \mc{F}^{\bullet}) = \sum_{S \in \mc{S}} \chi(S) \chi(\mc{H}^{\bullet}(\mc{F}^{\bullet})_{x_S}), \eeq where $x_S \to S$ in the inclusion of an arbitrary point in a stratum $S \in \mc{S}$ and $\bb{H}(-)$ is hypercohomology. In particular, letting $\mc{F}^{\bullet}= P^{\bullet}_{\tau}(K)$, letting $X=\mathscr{X}^{\tau}_{irr}(K)$ and fixing a stratification $\mc{S}$ with respect to which $P^{\bullet}_{\tau}(K)$ is constructible, we find that 
\begin{align} \label{equation:behrend-equality}
\chi_{\tau}(K):= \sum (-1)^i \bb{H}^i(X, P^{\bullet}_{\tau}(K)) &=\sum_{S \in \mc{S}} \chi(S) \chi(\mc{H}^{\bullet}(P^{\bullet}_{\tau}(K))_{x_S})  \\
&=\sum_{S \in \mc{S}} \chi(S) (-1)^{\op{dim} V}(1-\op{MF}_{f}(x_S)) \nonumber\\
&= \sum_{S \in \mc{S}} \chi(S) \nu_X(x_S)\nonumber \\
&= \chi(X, \nu_X) =\chi_B(X) =: \chi_B(\mathscr{X}^{\tau}_{irr}(K)) \nonumber
\end{align}

\section{Generic independence of the weight}

The purpose of this section is to prove \Cref{theorem:cassonconstant} from the introduction. In fact, most of the effort is directed at proving \Cref{theorem:independence}, which is purely a statement about algebraic geometry. We deduce \Cref{theorem:cassonconstant} as an easy corollary in \Cref{subsection:proof-thm-A}. 

\subsection{Setup and algebraic preliminaries}

Given a scheme $X$ over $\C$, recall that $\chi_B(X)$ is the Euler characteristic of $X$ weighted by the Behrend function. 

Let $X \subset \bb{A}^n$ be an affine $\C$-scheme corresponding to the ideal $I \subset \C[x_1,\dots,x_n]$. Let $X^0 \subset X$ be an open embedding. For $\tau \in \bb{A}^1$ a (closed) point, let $X^0_{\tau} \sub X_{\tau}$ be the (scheme-theoretic) fiber of the morphism to $\bb{A}^1$ induced by the projection $(x_1,\dots,x_n) \mapsto x_1$. 

\thm \label{theorem:independence} There is a (Zariski) open subset $\mc{V} \sub \bb{A}^1$ such that $\chi_B(X^0_{\tau})$ is constant for all $\tau \in \mc{V}$. \ethm

The proof of \Cref{theorem:independence} will occupy the next three sections.

Let us view $\C[x_1,\dots,x_n]/I$ as a $\C[x]$-algebra via the natural map taking $x \mapsto x_1$.  Observe that there are then isomorphisms 
\begin{align}
\lt( \C[x_1,\dots,x_n]/I \otimes_{\C[x]} \C[x]/(x-\tau) \rt) &= \C[x_1,\dots,x_n]/ (I, x_1-\tau) \label{equation:isom-0} \\
&= \C[z_2,\dots,z_n]/I_{\tau}, \nonumber
\end{align}
where the second map takes $(x_1,x_2,\dots,x_n) \mapsto (\tau, z_2,\dots,z_n)$ and $I_{\tau}$ is the image of $I$ in $\C[z_2,\dots,z_n]$. 

Hence we have that $$X_{\tau}= \op{Spec}(\C[z_2,\dots,z_n]/I_{\tau}) \sub \bb{A}^{n-1}.$$

Let us set $$R= \bigoplus_{m \geq 0} I^m/I^{m+1}$$ and observe that $R$ is naturally a $\C[x_1,\dots,x_n]/I$-algebra via the inclusion of the zero-graded piece. 

Let $$C= C_{X/ \bb{A}^n}= \op{Spec} R, $$ be the normal cone of $X \sub \bb{A}^n$ and consider the map 
\eq \label{equation:phi} \phi: \bigoplus_{m \geq 0} I^m/I^{m+1} \otimes_{\C[x]} \C[x]/(x-\tau) \to \bigoplus_{m \geq 0} I_{\tau}^m/I_{\tau}^{m+1}.\eeq

\prop \label{proposition:isom} Suppose that $x_1-\tau$ is not a zero-divisor in $\C[x_1,\dots,x_n]/I^l$ for any $l$. Then $\phi$ is an isomorphism. \eprop

\pf We have a natural isomorphism 
\eq  \bigoplus_{m \geq 0} I^m/I^{m+1} \otimes_{\C[x]} \C[x]/(x-\tau)  = \bigoplus_{m \geq 0} I^m/I^m(I+(x_1-\tau)), \eeq
which allows us to rewrite $\phi$ as the natural projection map
\eq \phi: \bigoplus_{m \geq 0} I^m/I^m(I+(x_1-\tau)) \to \bigoplus_{m \geq 0} I_{\tau}^m/I_{\tau}^{m+1}. \eeq

This map is clearly surjective. To check injectivity, choose $\a \in I^m$ and suppose that the composition $$I^m \to I^m/I^m(I+(x_1-\tau)) \to I_{\tau}^m/ I_{\tau}^{m+1}$$ annihilates $\a$. Then there exists $\b \in I^{m+1}$ such that $\a-\b \in (x_1-\tau)$. Hence $\a-\b \in (x_1-\tau) \cap I^m$.  The proposition now follows from \Cref{lemma:intersect} below. \epf

\lem \label{lemma:intersect}  Suppose that $x_1-\tau$ is not a zero-divisor in $\C[x_1,\dots,x_n]/I^l$ for any $l \geq 0$.  Then in $\C[x_1,\dots,x_n]$, we have $(x_1-\tau) \cap I^m= (x_1-\tau)I^m$ for all $m\geq 0$. \elem
\pf We only need to check the nontrivial inclusion. Any element $(x_1-\tau) \cap I^m$ is of the form $(x_1-\tau)\g$ for some $\g \in \C[x_1,\dots,x_n]$. But $(x_1-\tau) \g \in I^m$ implies that $(x_1-\tau) \g =0$ in $\C[x_1,\dots,x_n]/I^m$. By hypothesis, this implies $\g \in I^m$ and therefore $(x_1-\tau) \g \in (x_1-\tau) I^m$. \epf

Our next task is to check that the assumptions of \Cref{proposition:isom} and \Cref{lemma:intersect} are satisfied generically. 

\lem \label{lemma:gradeddivisor} Suppose that $f \in \C[x_1,\dots,x_n]$ is a zero-divisor of $\C[x_1,\dots,x_n]/I^l$ for some $l \geq 0$. Then $f$ (viewed as an element of the $0$-graded piece of $R$) is a zero-divisor in $R$.  \elem

\pf By hypothesis, there exists $g \in \C[x_1,\dots,x_n]$ such that $g \notin I^l$ but $fg \in I^l$. Let $0 \leq k <l$ be the largest integer such that $g \in I^k$ but $g \notin I^{k+1}$. Observe that $g$ can be viewed as a non-zero element of the $k$-graded piece of $R$. Viewing $f$ as an element of the $0$-graded piece of $R$, we have $0 = fg \in R$. Hence $f$ is a zero-divisor in $R$. \epf

\cor \label{corollary:zerodivisor} For all but finitely many $\tau \in \C$, the element $x_1-\tau \in \C[x_1,\dots,x_n]$ is not a zero divisor in the quotient ring $\C[x_1,\dots,x_n]/I^l$ for any $l \geq 0$.  \ecor

\pf Observe that $R$ is generated in degrees $0$ and $1$ as a $\C$-algebra, so $R$ is in particular a finitely generated $\C$-algebra.  In particular, $R$ is a Noetherian ring and it therefore has finitely many associated prime ideals whose union is precisely the set of zero-divisors of $R$.  If we suppose for contradiction that the corollary is false, then it follows from \Cref{lemma:gradeddivisor} that $(x_1-\tau)$ is a zero-divisor in $R$ for infinitely many values of $\tau \in \C$.  Hence there exist $\tau_1,\tau_2 \in \C$ with $\tau_1 \neq \tau_2$ such that $(x_1-\tau_1)$ and $(x_1-\tau_2)$ are both elements of the same associated prime ideal. Since $(x_1-\tau_1)-(x_1-\tau_2)=(\tau_2-\tau_1)$ is a unit, this gives the desired contradiction. \epf

\cor \label{corollary:coneisom} 
There is a Zariski open set $\mc{U}_1 \sub \bb{A}^1$ such that \eqref{equation:phi} is an isomorphism for all $\tau \in \mc{U}_1$. We therefore have the following commutative diagram of schemes, for all $\tau \in \mc{U}_1$: 
\eq \label{equation:diagram-0}
\begin{tikzcd}
C_{\tau} \arrow[r, "\sim"] \arrow{d} & C_{X_{\tau}/\bb{A}^{n-1}} \arrow{d} \\
X_{\tau} \arrow[r, "\sim"] & X_{\tau}.  
\end{tikzcd}
\eeq
\ecor

\pf The fact that \eqref{equation:phi} is an isomorphism for $\tau \in \mc{U}_1$ follows from \Cref{proposition:isom} and \Cref{corollary:zerodivisor}. The restriction of \eqref{equation:phi} to the zero-graded piece is just the isomorphism \eqref{equation:isom-0}, so we get the above diagram of schemes by taking $\op{Spec}(-)$.  \epf

\subsection{Passage to a cover} In general, the irreducible components of $C$ are not in bijection with the irreducible components of the fibers $C_{\tau}$ associated to the projection $C \to X \to \bb{A}^1$. However, the next proposition shows that this property becomes true in an open subset after passing to a suitable branched cover of $\bb{A}^1$, and that the fibers can moreover be assumed to have generically constant multiplicity and dimension. These facts are well-known in algebraic geometry, but we provide a detailed argument for completeness. 

\prop \label{proposition:cover} There exists an open set $\mc{U} \sub \bb{A}^1$ and a finite \'{e}tale cover $\psi: \tilde{\mc{U}} \to \mc{U}$ such that the following holds: if we let $Q^1,\dots,Q^q$ be the irreducible components of $C_{\tilde{\mc{U}}}$, then for all $p \in \tilde{\mc{U}}$, the irreducible components of $C_p$ are precisely $Q^1_p,\dots,Q^q_p$.  Moreover, the multiplicity and dimension of the $Q^i_p$ is independent of $p \in \tilde{\mc{U}}$.  \eprop 

The proof is an immediate consequence of the next three lemmas. Before stating these lemmas, it will be useful to make the following remark. 

\rmk  \label{remark:coveringshrink} If $p: \tilde{\mc{U}} \to \mc{U}$ is a finite \'{e}tale map between smooth $\C$-schemes of dimension $1$, then given an open subset $\mc{U}' \sub \mc{U}$, the map $p$ restricts to an \'{e}tale map $p|_{\mc{U}'} : p^{-1}(\mc{U}') \to \mc{U}'$. If $\mc{V} \sub \tilde{\mc{U}}$ is an open subset, then after taking a possibly smaller open $\mc{V}' \sub \mc{V} \sub \tilde{\mc{U}}$, we can assume that $p$ restricts to an \'{e}tale map onto its image. Indeed, observe that $\tilde{\mc{U}}- \mc{V}$ is closed. Hence $p(\tilde{\mc{U}}- \mc{V})$ is closed (since finite morphisms are closed) and $p^{-1}(p(\tilde{\mc{U}}- \mc{V}))$ is also closed.  We can then set $\mc{V}' = \tilde{\mc{U}} - p^{-1}(p(\tilde{\mc{U}}- \mc{V}))$. \ermk


\lem \label{lemma:step1} There exists an open set $\mc{U} \sub \bb{A}^1$ and a finite \'{e}tale cover $\psi: \tilde{\mc{U}} \to \mc{U}$ such that the irreducible components of the generic fiber of $C_{\tilde{\mc{U}}} \to \tilde{\mc{U}}$ are geometrically irreducible.  \elem

\pf Let $\eta \in \bb{A}^1$ be the generic point and note that $k(\eta)\simeq \C(x)$. By applying \cite[\href{https://stacks.math.columbia.edu/tag/054R}{Tag 054R}]{stacks-project}, there is a finite extension $K/ k(\eta)$ such that $C_K$ is geometrically irreducible over $K$. Observe that $K$ has transcendence degree one over $\C$. Hence, according to \cite[\href{https://stacks.math.columbia.edu/tag/0BY1}{Tag 0BY1}]{stacks-project}, this extension is induced by a dominant rational map $f: \S \to \bb{A}^1$, where $\S$ is an algebraic curve over $\C$. By the theorem on generic smoothness on the target \cite[25.3.3]{vakil} (and the fact that $f$ is dominant), there is an open $\mc{U} \sub \bb{A}^1$ such that $f|_{f^{-1}(\mc{U})}$ is smooth. Since $f$ is evidently of relative dimension zero, $f|_{f^{-1}(\mc{U})}$ is \'{e}tale. The lemma follows with $\tilde{\mc{U}}:= f^{-1}(\mc{U})$ and $\psi:=f$. \epf

\lem \label{lemma:step2} Let $C_{\tilde{\mc{U}}} \to \tilde{\mc{U}}$ be as in \Cref{lemma:step1}. After possibly shrinking $\tilde{\mc{U}}$ (cf. \Cref{remark:coveringshrink}), we can assume that the following holds: if we let $Q^1,\dots,Q^q$ be the irreducible components of $C_{\tilde{\mc{U}}}$, then for all $p \in \tilde{\mc{U}}$, the irreducible components of $C_p$ are precisely $Q^1_p,\dots,Q^q_p$. \elem

\pf Let $C^1_K,\dots,C^q_K$ be the irreducible components of the generic fiber $C_K$. For $1 \leq i \leq q$, let $C^i$ be the smallest closed irreducible subscheme of $C_{\tilde{\mc{U}}}$ whose generic fiber is $C^i_K$.  According to \cite[\href{https://stacks.math.columbia.edu/tag/054Y}{Tag 054Y}]{stacks-project}, we can assume after possibly shrinking $\tilde{\mc{U}}$ that the irreducible components of the fiber $C_p$ for any $p \in \tilde{\mc{U}}$ are precisely $\{C_p^{i,j}\}_{j=1}^{n_i}$ where $1 \leq i \leq q$. It now follows from \Cref{lemma:step1} and \cite[\href{https://stacks.math.columbia.edu/tag/0559}{Tag 0559}]{stacks-project} (geometric irreducibility spreads out) that $n_i=1$ for all $i$. This completes the proof.  (Note that geometric irreducibility plays a crucial role, since the analog of \cite[\href{https://stacks.math.columbia.edu/tag/0559}{Tag 0559}]{stacks-project} is false for irreducible schemes which are not geometrically irreducible). \epf

\lem \label{lemma:step3} Let $C_{\tilde{\mc{U}}} \to \tilde{\mc{U}}$ satisfy the conditions of \Cref{lemma:step1} and \Cref{lemma:step2}. After possibly shrinking $\tilde{\mc{U}}$, we can assume that the multiplicity and dimension of the $Q^i_p$ is independent of $p \in \tilde{\mc{U}}$.  \elem
\pf \Cref{lemma:step2} gives a bijection between the irreducible components of the generic fiber $C_K$ and the irreducible components of the fibers $C_p$ for $p \in \tilde{\mc{U}}$. The present lemma is simply a consequence of the fact that both dimension and multiplicity ``spread out"; that is, after possibly further shrinking $\tilde{\mc{U}}$, we can assume that the bijection constructed in \Cref{lemma:step2} preserves dimension and multiplicity. The relevant reference for dimension is \cite[\href{https://stacks.math.columbia.edu/tag/02FZ}{Tag 02FZ}]{stacks-project}; for multiplicity, one can apply \cite[III, 9.8.6]{ega4} to the structure sheaf of $C$. Note that the notion of geometric multiplicity in \cite[III, 9.8.6]{ega4} agrees with our notion of multiplicity since we are in characteristic zero; see \cite[II, 4.7.5]{ega4}. \epf

\pf[Proof of \Cref{proposition:cover}] Combine \Cref{lemma:step2} and \Cref{lemma:step3}.\epf

Note that the image of each $Q^i$ under the map $C_{\tilde{\mc{U}}} \to X_{\tilde{\mc{U}}}$ is irreducible (since the image of an irreducible set under a continuous map is irreducible). It is also closed: this follows by combining \cite[B.5.3.]{fulton} and the fact that $\psi: \tilde{\mc{U}} \to \mc{U}$ is an \'{e}tale cover. We let $V^i$ be the image of $C^i$ and conclude that $V^i$ is an irreducible subvariety of $X_{\tilde{\mc{U}}}$ when endowed with the canonical reduced closed subscheme structure. Observe also that the fibers $V^i_p$ are irreducible for $p \in \tilde{\mc{U}}$. Indeed, the $Q^i_p$ are irreducible, so this follows from the fact that $V^i_p= \pi(Q^i_p)$.

\subsection{Stratification theory}   All stratifications which we consider in this section will be assumed to be w-stratifications in the sense of \Cref{definition:wstrat}. In particular, this implies that our stratifications are Whitney stratifications and that the strata are smooth, connected, locally closed subvarieties.

\defi[cf. (3.2) in \cite{verdier}] \label{definition:transverse} Given a morphism $f: X \to Y$ of complex algebraic varieties and a stratification $\mc{S}$ of $X$, we say that $f$ is transverse to $\mc{S}$ if $f$ restricts to a smooth morphism on each stratum. \edefi

As observed in \cite[(3.6)]{verdier}, if $f:X \to Y$ is transverse to a stratification $\mc{S}$, then given any $y \in Y$, the fiber $f^{-1}(y)$ inherits a stratification by restriction of the strata. 

It will be convenient to record the following lemma, whose proof is a routine verification. 

\lem \label{lemma:restrictstrat} Suppose that $V' \sub V$ is a (locally closed) subvariety of $V$. Suppose that $\mc{S}$ is a stratification of $V$ such that $V'$ is a union of strata. Then $\mc{S}|_{V'}$ is a stratification of $V'$ (in particular, $\mc{S}$ also satisfies the axioms of \Cref{definition:wstrat}). \elem

We consider $\tilde{\pi}: X_{\tilde{\mc{U}}} \to \tilde{\mc{U}}$ satisfying the properties of \Cref{proposition:cover}.   

\prop \label{proposition:product} After possibly replacing $\mc{U}$ with a smaller open $\mc{U}_2 \sub \mc{U}$ (cf. \Cref{remark:coveringshrink}), we can assume that $X_{\tilde{\mc{U}}}$ admits a stratification $\mc{S}$ with the following properties:
\begin{itemize}
\item[(i)] The stratification $\mc{S}$ is transverse to $\tilde{\pi}$. 
\item[(ii)] The subvarieties $X^0_{\tilde{\mc{U}}}, V^1,\dots,V^q$ of $X_{\tilde{\mc{U}}}$ are a union of strata. 
\item[(iii)] For each $x \in \tilde{\mc{U}}$, there exists a ball $B_x \sub \tilde{\mc{U}}$ such that $\tilde{\pi}^{-1}(B_x)$ is homeomorphic to $X_x \tms B_x$. Moreover, this homeomorphism is compatible with the projection and preserves the natural product stratification. 
\end{itemize}
\eprop

\pf  For ease of notation, we write $f= \tilde{\pi}: X_{\tilde{\mc{U}}} \to \tilde{\mc{U}}$. We will argue exactly as in the proof of Proposition 5.1 in \cite{verdier}.  Applying the Nagata compactification theorem, we can factor $f: X_{\tilde{\mc{U}}} \to \tilde{\mc{U}}$ as an open embedding $i: X_{\tilde{\mc{U}}} \to \ov{X_{\tilde{\mc{U}}}}$ followed by a proper map $\ov{f}: \ov{X_{\tilde{\mc{U}}}} \to \tilde{\mc{U}}$.  Given an open set $V \sub \tilde{\mc{U}}$, we write $f|_{V}$ or $\ov{f}|_{V}$ for the restriction of $f$ or $\ov{f}$ to $f^{-1}(V)$ or $\ov{f}^{-1}(V)$ respectively.  

According to (2.2) in \cite{verdier}, we can choose a Whitney stratification $\mc{S}$ of $\ov{X_{\tilde{\mc{U}}}}$ so that $X^0_{\tilde{\mc{U}}}, V^1,\dots,V^q$ and $X_{\tilde{\mc{U}}}$ are a union of strata.  Next, (3.3) in \cite{verdier} shows that one can find an open $\mc{V} \sub \tilde{\mc{U}}$ so that $\ov{f}$ is transverse on $\ov{f}^{-1}(\mc{V})$ to $\mc{S} \cap \ov{f}^{-1}(\mc{V})$.   Finally, Verdier shows in (4.14) of \cite{verdier} that there are trivializations of $\ov{f}|_{\mc{V}}$ with the desired properties, i.e. which are compatible with projection and preserve the stratifications.  Since $X_{\tilde{\mc{U}}}$ is a union of strata, it follows that these also give local trivializations for $f|_{\mc{V}}$, as desired. \epf

Since $X^0_{\tilde{\mc{U}}}, V^1,\dots,V^s$ are each a union of strata (and since $\mathring{V}^i := X^0_{\tilde{\mc{U}}} \cap V^i$ is therefore also a union of strata), we obtain the following corollary of (ii) and (iii): 

\cor \label{corollary:product} If we let $\tilde{\pi}^0 :X^0_{\tilde{\mc{U}}} \to \tilde{\mc{U}}$ be the composition $X^0_{\tilde{\mc{U}}} \hookrightarrow X_{\tilde{\mc{U}}} \to \tilde{\mc{U}}$, then (iii) holds for $X^0_{\tilde{\mc{U}}}, \tilde{\pi}^0, W^0_{x}$ in place of $X_{\tilde{\mc{U}}}, \tilde{\pi}, X_x$. If, for $i=1,\dots,v$, we let $\tilde{\pi}^i: V^i \to \mc{U}$ be the composition $V^i \hookrightarrow X_{\tilde{\mc{U}}} \to \mc{U}$, then (iii) holds for $V^i, \tilde{\pi}^i, V^i_x$ in place of $X_{\tilde{\mc{U}}}, \tilde{\pi}, X_x$. Finally, letting $\mathring{V}^i := X^0_{\tilde{\mc{U}}} \cap V^i$ and letting $\tilde{\pi}^{i,0}: \mathring{V}^i \hookrightarrow X_{\tilde{\mc{U}}} \to \mc{U}$ be the obvious composition, then (iii) holds for $\mathring{V}^i, \tilde{\pi}^{i,0}, \mathring{V}^i_x$ in place of $X_{\tilde{\mc{U}}}, \tilde{\pi}, X_x$.  \ecor

\subsection{Completion of the argument} We now have the ingredients in place to prove \Cref{theorem:independence}. For $\tau \in \bb{A}^1$, it follows from \cite[Prop.\ 1.5(i)]{behrend} that $\nu_{X_{\tau}^0}= \nu_{X_{\tau}}|_{X_{\tau}^0}$. Hence, we have
\eq \chi_B(X_{\tau}^0):= \chi(X_{\tau}^0, \nu_{X_{\tau}^0}) = \chi(X_{\tau}^0, \nu_{X_{\tau}}|_{X_{\tau}^0}). \eeq

According to \Cref{corollary:coneisom} and \Cref{proposition:cover}, up to replacing $\mc{U}_1$ and $\mc{U}$ by a possibly smaller open set $\mc{U}_3 \sub \mc{U}_1 \cap \mc{U}$, we can assume that for all $\tau \in \mc{U}$ and $p \in \tilde{\mc{U}}$ satisfying $\psi(p)= \tau$ there is a diagram
\eq \label{equation:diagram}
\begin{tikzcd}
C_p \arrow[r, "\sim"] \arrow[d] & C_{\tau} \arrow[r, "\sim"] \arrow{d} & C_{X_{\tau}/\bb{A}^{n-1}} \arrow{d} \\
X_p \arrow[r, "\sim"] & X_{\tau} \arrow[r, "\sim"] & X_{\tau}.  
\end{tikzcd}
\eeq
 
We then have, by \eqref{equation:canonicalcycle} and \Cref{proposition:cover},
$$\nu_{X_{p}} = \op{Eu} \lt( \sum_{i=1}^q a_i(p) V_p^i \rt), $$
where \eq \label{equation:formulaai} a_i(p)= (-1)^{\op{dim} \pi(Q^i_p)} \op{mult}(Q^i_p) = (-1)^{\op{dim} V^i_p} \op{mult}(Q^i_p). \eeq

We conclude that 
\eq \label{equation:tautop} \chi_B(X^0_{\tau}) = \chi \lt(X_{p}^0, \op{Eu} ( \sum_{i=1}^q a_i(p) V^i_p)\big|_{X_{p}^0} \rt). \eeq

By appealing to the complex-analytic definition of the local Euler obstruction (see \cite[Sec.3]{macpherson} or \cite[p.\ 32]{joyce-song}), we see that the local Euler obstruction of a cycle at some point $x$ only depends on an analytic neighborhood of $x$. Hence we have
\begin{align}
\op{Eu}(\sum a_i(p) V_p^i)|_{X_p^0} &= \op{Eu}( \sum a_i(p) ( V^i_p \cap X_p^0)) \nonumber \\
&= \sum_{\{j \in \S\}} a_j(p) \op{Eu}(\mathring{V}^j_p) \label{equation:euler-local},
\end{align}
where $\S \sub \{1,2\dots,q\}$ and $j \in \S$ iff $V^j \cap X_{\tilde{\mc{U}}}^0 = \emptyset$, and where we let $\mathring{V}_p^j = V_p^i \cap X_p^0$. Note that the second line follows from the fact that $\op{Eu}(-)$ is a homomorphism from the group of algebraic cycles on $X_p^0$ to the group of constructible functions; see \cite[p.\ 376]{fulton}.

Next, it follows from \cite[1.3(ii)]{behrend} that 
\eq \label{equation:b} \chi(X_p^0, \sum_{j \in \S} a_j(p) \op{Eu}(\mathring{V}^j_p)) = \sum_{j \in \S} a_j(p) \chi(X_p^0, \op{Eu}(\mathring{V}^j_p))= \sum_{j \in \S} a_j(p) \chi(\mathring{V}^j_p, \op{Eu}(\mathring{V}^j_p)). \eeq

\prop \label{proposition:localconstcoeff} 
	After possibly replacing $\tilde{\mc{U}}$ with a smaller open subset $\mc{U}_4$, we can assume that the function $p \mapsto a_i(p)$ is constant for $p \in \tilde{\mc{U}}$ for all $i=1,\dots,q$. \eprop
\pf Noting that $V^i_p=\pi(Q^i_p)= \pi(Q^i)_p$, it follows from \cite[\href{https://stacks.math.columbia.edu/tag/05F7}{Tag 05F7}]{stacks-project} that $\op{dim} V^i_p$ is constant on an open subset of $\tilde{\mc{U}}$. The result now follows by combining \eqref{equation:formulaai} and \Cref{proposition:cover}. \epf

We will also need the following lemma: 
\lem[Lem.\ 1.1(3) in \cite{paru-pra}] \label{lemma:paru-pra} Assume that an irreducible variety $Y$ is embedded in $\C^n$ and a nonsingular subvariety $Z$ intersects a Whitney stratification of $Y$ transversally. Then $\op{Eu}(Z \cap Y)(x)= \op{Eu}(Y)(x)$ for all $x \in Z \cap Y$. \elem

\prop \label{proposition:localconstchi} After possibly replacing $\tilde{\mc{U}}$ with a smaller open $\mc{U}_5 \sub \tilde{\mc{U}}$, we may assume that the function $p \mapsto \chi(\mathring{V}^i_p, \op{Eu}(\mathring{V}^i_p))$ is locally (and hence globally) constant for $p \in \tilde{\mc{U}}$ and $i=1,\dots,v$. \eprop
\pf
It follows from \Cref{proposition:product}(ii) that $V^i$ is a union of strata of the stratification $\mc{S}$ of $X_{\tilde{\mc{U}}}$.  \Cref{lemma:restrictstrat} then implies that the restriction of $\mc{S}$ to $\mathring{V}^i$ is a stratification which we call $\mathring{\mc{S}}^i$. According to \Cref{proposition:product}(i) and the comment following \Cref{definition:transverse}, the fiber $\mathring{V}^i_p$ inherits a stratification from $\mathring{\mc{S}}^i$ which we call $\mathring{\mc{S}}_p^i$. Note that the strata of $\mathring{\mc{S}}_p^i$ are of the form $S_p$ for $S \in \mathring{\mc{S}}^i$. 

According to \cite[Prop.\ 2]{brasselet}, the constructible function $\op{Eu}(\mathring{V}^i)$ is constant on each stratum of $\mathring{\mc{S}}^i$. It follows that $\op{Eu}(\mathring{V}_p^i)|_{S'}$ is constant for any $S' \in \mathring{\mc{S}}_p^i$. If $S'= S_p$ for some stratum $S \in \mathring{\mc{S}}^i$, then \Cref{proposition:product}(i) and \Cref{lemma:paru-pra} imply that $\op{Eu}(\mathring{V}_p^i)|_{S'} = \op{Eu}(\mathring{V}^i)|_S$ is constant. 

It now follows that
\begin{align} 
\chi(\mathring{V}^i_p, \op{Eu}(\mathring{V}^i_p)) &:= \sum_{S' \in \mathring{\mc{S}}_p^i} \chi(S') \op{Eu}(\mathring{V}^i_p)|_{S'} \nonumber\\
&=\sum_{S \in \mathring{\mc{S}}^i} \chi(S_p) \op{Eu}(\mathring{V}^i)|_{S} \label{equation:sum-indep}. 
\end{align}

According to \Cref{corollary:product}, the topological Euler characteristic of the fiber $\chi(S_p)$ is independent of $p \in \tilde{\mc{U}}$, for $S \in \mathring{\mc{S}}^i$. It follows that \eqref{equation:sum-indep} is independent of $p \in \tilde{\mc{U}}$, which is what we wanted to show.  \epf 

\pf[Proof of \Cref{theorem:independence}]   By combining \eqref{equation:tautop}, \eqref{equation:euler-local} and \eqref{equation:b}, we find that \eq \chi_B(X^0_{\tau})  =  \sum_{j \in \S} a_j(p) \chi(\mathring{V}^j_p, \op{Eu}(\mathring{V}^j_p)). \eeq Setting $\mc{V}= \tilde{\mc{U}}$, it follows from \Cref{proposition:localconstcoeff} and \Cref{proposition:localconstchi} that this expression is constant for $p \in \mc{V}$. \epf

\subsection{Proof of \Cref{theorem:cassonconstant}} \label{subsection:proof-thm-A}
Fix a presentation $\G = \langle g_1,\dots,g_m \mid r_1,\dots,r_l\rangle$. As explained for example in \cite[Sec.\ 2.1]{cote-man}, one can associate to this presentation a $\C$-algebra $\mc{A}(\G)$ with generators $x_{11}^{g_1}, x_{12}^{g_1}, x_{21}^{g_1}, x_{22}^{g_1}, \dots, x_{11}^{g_m}, x_{12}^{g_m}, x_{21}^{g_m}, x_{22}^{g_m}$ whose spectrum is the representation scheme $\mathscr{R}(\G)$.  The group scheme $\op{SL}_2$ acts on $\mathscr{R}(\G)$ by conjugation. Letting $R$ be the ring of functions of $\op{SL}_2$, this action is induced by a map $\mu: \mc{A}(\G) \to \mc{A}(\G) \otimes R$. The ring of invariants is $\mc{A}(\G)^{\op{SL}_2}$ is the character scheme is $\op{Spec} \mc{A}(\G)^{\op{SL}_2}$. 

Note that $\mc{A}(\G)^{\op{SL}_2}$ is a finitely-generated $\C$-algebra, so we can fix a generating set $X_1,\dots,X_n$ and we can moreover assume that $X_1= x^{g_1}_{11}+x^{g_1}_{22}$. The surjective ring map $\C[x_1,\dots,x_n] \to \mc{A}(\G)^{\op{SL}_2}$ sending $x_i \mapsto X_i$ induces an isomorphism 
\eq \C[x_1,\dots,x_n]/I \to \mc{A}(\G)^{\op{SL}_2},\eeq where $I$ is the kernel of the surjection.  This gives an embedding of schemes 
$$\mathscr{X}(\G) \hookrightarrow \bb{A}^n.$$ There is also an open embedding $\mathscr{X}_{irr}(\G) \hookrightarrow \mathscr{X}(\G)$; see \cite[p.\ 7]{abou-man}. After composing with the projection $\bb{A}^n \to \bb{A}^1$ sending $(x_1,\dots,x_n) \mapsto x_1$, we get a morphism $\mathscr{X}(\G) \to \bb{A}^1$. 

According to \Cref{theorem:independence} applied to $\mathscr{X}_{irr}(\G) \sub \mathscr{X}(\G) \sub \bb{A}^n$, there is a Zariski open $\mc{V} \sub \bb{A}^1$ such that $\chi_B(\mathscr{X}_{irr}(\G)_{\tau})$ is constant over all $\tau \in \mc{V}$. To complete the proof of \Cref{theorem:cassonconstant}, it follows from \eqref{equation:behrend-equality} that it is enough to prove that $\mathscr{X}_{irr}(\G)_{\tau}= \mathscr{X}^{\tau}_{irr}(\G)$ for all but finitely many $\tau \in \bb{A}^1$. This is the content of the following proposition.

\prop \label{proposition:fiberisom} \begin{itemize}
\item[(i)]$\mathscr{X}(\G)_{\tau}= \mathscr{X}^{\tau}(\G)$ for all but finitely many values of $\tau \in \bb{A}^1$. 
\item[(ii)] $\mathscr{X}_{irr}(\G)_{\tau}= \mathscr{X}^{\tau}_{irr}(\G)$ for all but finitely many values of $\tau \in \bb{A}^1$.
\end{itemize}
 \eprop

\pf[Proof of (i)]  For $\tau \in \bb{A}^1$, let $h_t = x^{g_1}_{11}+ x^{g_1}_{22} - t$. Thus $h_{\tau}$ is the image of $x_1-\tau$ under the isomorphism $\C[x_1,\dots,x_n]/I \to \mc{A}(\G)^{\op{SL}_2}$ described above. Consider the surjective map 
\eqs \mc{A}(\G)^{\op{SL}_2} \to ( \mc{A}(\G) / h_{\tau} \mc{A}(\G))^{\op{SL}_2} = \mc{A}(\G)^{\op{SL}_2} / ( h_{\tau} \mc{A}(\G))^{\op{SL}_2},\eeqs  
which induces the quotient map
\eq \label{equation:invariantsurj} \mc{A}(\G)^{\op{SL}_2} / (h_{\tau} \mc{A}(\G)^{\op{SL}_2}) \to \mc{A}(\G)^{\op{SL}_2} / (h_{\tau} \mc{A}(\G))^{\op{SL}_2}. \eeq

We wish to show that the morphism in \eqref{equation:invariantsurj} is injective for all but finitely many $\tau \in \bb{A}^1$. We closely follow the proof of \cite[Prop.\ 5.3]{cote-man}. To this end, observe that it is enough to establish the following containment: 
\eq \label{equation:containment} (h_{\tau} \mc{A}(\G))^{\op{SL}_2} \sub h_{\tau} \mc{A}(\G)^{\op{SL}_2}. \eeq

Let $R$ be the coordinate ring of the group scheme $\op{SL}_2$. Let $\mu: \mc{A}(\G) \to \mc{A}(\G) \otimes R$ be the $\C$-algebra morphism inducing the $\op{SL}_2$-action on $\mathscr{R}(\G)= \op{Spec}\mc{A}(\G)$.  By definition, $f \in \mc{A}(\G)^{\op{SL}_2}$ if and only if $\mu(f)= f \otimes 1$. 

Suppose for contradiction that \eqref{equation:containment} is false for infinitely many values of $\tau$. Then there exists $g \in \mc{A}(\G)$ such that $h_{\tau} g \in (h_{\tau} \mc{A}(\G))^{\op{SL}_2}$ but $g \notin \mc{A}(\G)^{\op{SL}_2}$. Hence 
\eq 0 = \mu(h_{\tau} g)- \mu(h_{\tau})\mu(g) = (h_{\tau} g \otimes 1) - (h_{\tau} \otimes 1) \mu(g)= (h_{\tau} \otimes 1)(g\otimes 1 - \mu(g)). \eeq  Since $g \notin \mc{A}(\G)^{\op{SL}_2}$, we have that $g \otimes 1 - \mu(g) \neq 0$ which implies that $(h_{\tau} \otimes 1)$ is a zero-divisor in the ring $\mc{A}(\G) \otimes R$. 

Since $\mc{A}(\G) \otimes R$ is Noetherian, it has finitely many associated prime ideals. Moreover, it is a general fact that every zero-divisor must be contained in one of these ideals; see \cite[(5.5.10)]{vakil}. By combining this fact with the previous paragraph, it follows that we can find $s, s' \in \bb{A}^1$ with $s \neq s'$ such that $(h_s \otimes 1)$ and $(h_{s'} \otimes 1)$ are contained in the same associated prime ideal. However, observe that $(h_s \otimes 1)- (h_{s'} \otimes 1) = (h_s-h_{s'} \otimes 1)= (s'-s \otimes 1)= (s'-s) (1 \otimes 1)$. This is a contradiction since $(s'-s) (1 \otimes 1)$ is a unit. \epf

\pf[Proof of (ii)]  The argument is similar to the proof of \cite[Prop.\ 5.5]{cote-man}. It follows from (i) that, for all but finitely many values of $\tau \in \bb{A}^1$, the composition $\mathscr{X}^{\tau}_{irr}(\G) \hookrightarrow \mathscr{X}^{\tau}(\G) = \mathscr{X}(\G)_{\tau}$ is an open embedding. On the other hand, $\mathscr{X}_{irr}(\G)_{\tau} \sub \mathscr{X}(\G)_{\tau}$ is also an open embedding.  It's clear that both open embeddings have the same closed points (corresponding to irreducible representations $\rho$ such that $\op{Tr}(\rho(g))=\tau$). Hence the claim follows from the fact that any two open subschemes which have the same closed points coincide. \epf  

\section{Additivity of (sheaf-theoretic) $\op{SL}(2,\C)$ Casson-Lin invariant} 

The goal of this section is to prove \Cref{theorem:connectedsum}, which states that the $\op{SL}(2,\C)$ Casson-Lin invariant $\chi_{CL}(K)$ is additive under connected sums of knots in integral homology $3$-spheres.

\subsection{The structure of the character variety} \label{subsection:structure-character-variety}

Following \cite[p.\ 5]{abou-man}, we will partition the $\op{SL}(2,\C)$ representations of a finitely-presented group $\G$ into five classes.  Let $B \sub \op{SL}(2,\C)$ be the Borel subgroup of upper-triangular matrices and let $B_P \sub B$ be the subgroup of matrices of the form $\pm \begin{pmatrix} 1 &a \\
0 & 1 
\end{pmatrix} $ for $a \in \C$. Let $D \sub \op{SL}(2,\C)$ be the group of diagonal matrices. Every representation $\G \to \op{SL}(2,\C)$ is of exactly one of the following types:

\begin{itemize}
\item[(a)] The irreducible representations. These representations have trivial stabilizer. 
\item[(b)] Representations which are conjugate to one in $B$ but not in $B_P$ or in $D$. 
\item[(c)] Representations which are conjugate to one in $B_P$ but not in $\{ \pm \op{Id}\}$. 
\item[(d)] Representations which are conjugate to one in $D$ but not in $\{ \pm \op{Id}\}$. 
\item[(e)] Representations with image in $\{ \pm \op{Id}\}$.
\end{itemize}

We say that representations of types (b)-(e) are \emph{reducible} and that representations of types (d) and (e) are \emph{abelian}. 

Let us now specialize to the case where $\G= \pi_1(K)$ for $K \sub Y$ a knot in an integral homology $3$-sphere. We fix $\tau \in \C-\{\pm 2\}$ and consider the relative character variety $X^{\tau}(K)$. The points of $X^{\tau}(K)$ correspond to irreducible representations $\rho: \G \to \op{SL}(2,\C)$ with $\op{Tr}(\rho(m))=\tau$ for $m \in \G$ a meridian.

It will be useful to introduce the following terminology. 
\defi
	Let $K \sub Y$ be as above. Let $\mc{G}(K) \sub \C$ be the set of values $\tau \in \C- \{ \pm 2\}$ with the property that $\tau \neq e^{\a/2} + e^{-\a/2}$ whenever $e^{\a}$ is a root of the Alexander polynomial of $K$, where $\a \in \C$. We say that representation $\rho: \pi_1(Y-K) \to \op{SL}(2,\C)$ is \emph{good} if $\op{Tr} \rho(m) \in \mc{G}(K)$ for $m \in \pi_1(Y-K)$ a meridian. 
\edefi

It is a remarkable fact first observed by de Rham (see \cite[Sec.\ 6.1]{planecurves}) that a good representation is reducible if and only if it is abelian.

\lem \label{lemma:disjoint-union}
	Let $K \sub Y$ be as above and suppose that $\tau \in \mc{G}(K) \sub \C- \{ \pm 2\}$. Then there is a (scheme-theoretic) decomposition 
$\mathscr{R}^{\tau}(K) = \mathscr{R}^{\tau}_{ab}(K) \sqcup \mathscr{R}^{\tau}_{irr}(K)$ where $\mathscr{R}^{\tau}_{ab}(K)$ is the image of $\mathscr{R}^{\tau}(\Z)$ under the abelianization map $\pi_1(Y-K) \to H_1(Y-K; \Z) \simeq \Z$. 
\elem

\pf
	The surjective map $\pi_1(Y-K) \to H_1(Y-K; \Z) \simeq \Z$ induces a closed embedding of relative representation schemes $\phi^{\tau}: \mathscr{R}^{\tau}(\Z) \hookrightarrow \mathscr{R}^{\tau}(K)$. Let $C^{\tau} \sub \mathscr{R}^{\tau}(K)$ be the union of all irreducible components of $\mathscr{R}^{\tau}(K)$ which contain a non-abelian representation. Since $\op{im} \phi^{\tau}$ and $C^{\tau}$ are closed and cover $\mathscr{R}^{\tau}(K)$, it's enough to show that they have empty intersection. 
	
	Suppose for contradiction that there is an abelian representation $\rho \in C^{\tau}$. Then $\rho$ belongs to an irreducible component $C_0^{\tau} \sub C^{\tau}$ which contains an non-abelian representation. Since all representations in $\mathscr{R}^{\tau}(K)$ are good, $C_0^{\tau}$ contains an irreducible representation $\rho'$. Let $\pi: \mathscr{R}^{\tau}(K) \to \mathscr{X}^{\tau}(K)$ be the natural projection map and note that the closure $\ov{\pi(C_0^{\tau})}$ is irreducible. It follows that there is an irreducible component $D_0 \sub \mathscr{X}(K)$ which contains $\ov{\pi(C_0^{\tau})}$, and hence contains both $\rho$ and $\rho'$. But this impossible in view of \cite[Prop.\ 6.2]{planecurves}, which states that if an abelian representation lies on a component of the character variety which also contains an irreducible representation, then this abelian representation is bad.  
	
	We have shown that $\mathscr{R}^{\tau}(K)= \mathscr{R}^{\tau}_{ab}(K) \sqcup C^{\tau}$, since both components are open and closed. It's clear that $C^{\tau}= \mathscr{R}^{\tau}_{irr}(K)$ since both open subschemes have the same closed points. Finally, the fact that the closed embedding $\mathscr{R}^{\tau}(\Z) \hookrightarrow \mathscr{R}^{\tau}_{ab}$ is an isomorphism (i.e. also an open embedding) is straightforward; cf. \cite[Prop.\ 7.6]{cote-man}.  \epf

The connected sum operation for knots is described in detail in \cite[Sec.\ 7.1]{cote-man}, and it will be useful to review this description in order to set our notation. 

Let $K_1 \sub Y_1$ and $K_2 \sub Y_2$ be oriented knots in integral homology $3$-spheres. Let $\ov{B}_1 \sub Y_1$ and $\ov{B}_2 \sub Y_2$ be small closed balls with the property that $\ov{B}_i - K_i$ is diffeomorphic to $\{ (x,y,z) \mid \|(x,y,z)\| \leq 1, (x,y,z) \neq (x,0,0) \}$ for $i=1,2$. Let $B_i \sub \ov{B}_i$ be the (open) interior. Let $C_i:= \d(Y_i-B_i)$ and let $\phi: C_1 \to C_2$ be an orientation reversing diffeomorphism which sends $\{ C_1 \cap K_1\} \to \{ C_2 \cap K_2 \}$ and preserves the orientation on the sets $\{C_i \cap K_i\}$ induced by the orientation of $K_i$. 

\defi 
	Let $Y=Y_1\# Y_2:= Y_1 \cup_{\phi} Y_2$ be obtained by gluing $Y_1-B_1$ to $Y_2-B_2$ via $\phi$ and let $K= K_1 \# K_2 := (K_1 - B_i) \cup_{\phi} (K_2- B_i)$ be the induced knot. We say that $K \sub Y$ is the \emph{connected sum} of $K_1$ and $K_2$. While this construction appears to depend on choices, it can be shown that $K \sub Y$ is well-defined up to equivalence of knots.  
\edefi

For the remainder of this section, we assume that $K_i \sub Y_i$ are fixed and let $K=K_1\#K_2$. By van Kampen's theorem, we have 
\begin{align} \label{equation:fiberconnect} \pi_1(Y-K) &= \pi_1(Y_1-K_1-B_1) *_{\pi_1(S^2-p_1-p_2)} \pi_1(Y_2-K_2-B_2) \\
&= \pi_1(Y_1-K_1) *_{\pi_1(S^2-p_1-p_2)} \pi_1(Y_2-K_2). \nonumber
\end{align}

Here, we have identified $\d (Y_1-B_1)= \d (Y_2-B_2)= (S^2-p-q)$ via $\phi$, for $p,q$ a pair of distinct points on $S^2$. We can assume that the above fundamental groups are computed with respect to some reference basepoint $x \in S^2-p-q$. 

Since the class of the meridian generates $\pi_1(S^2-p-q)$, we find that the representations of $\pi_1(Y-K)$ are pairs of representations $(\rho_1,\rho_2)\in \mathscr{R}^{\tau}(\pi_1(Y_1-K_1))\times\mathscr{R}^{\tau}(\pi_1(Y_2-K_2))$ such that $\rho_1$ and $\rho_2$ agree on the meridian. In fact, by combining \Cref{lemma:disjoint-union} and \eqref{equation:fiberconnect}, we get the following fiber product presentation for the relative representation scheme:
\begin{align}
\mathscr{R}^{\tau}(K) &= \lt( \mathscr{R}^{\tau}_{ab}(K_1) \sqcup \mathscr{R}^{\tau}_{irr}(K_1) \rt) \tms_{\mathscr{R}^{\tau}(\Z)}  \lt( \mathscr{R}^{\tau}_{ab}(K_2) \sqcup \mathscr{R}^{\tau}_{irr}(K_2)  \rt)  \nonumber \\ 
&= \lt( \mathscr{R}^{\tau}_{ab}(K_1) \tms_{\mathscr{R}^{\tau}(\Z)}  \mathscr{R}^{\tau}_{ab}(K_2)   \rt)  \sqcup  \lt( \mathscr{R}^{\tau}_{irr}(K_1) \tms_{\mathscr{R}^{\tau}(\Z)} \mathscr{R}^{\tau}_{ab}(K_2) \rt) \label{equation:decomp} \\& \sqcup  \lt( \mathscr{R}^{\tau}_{ab}(K_1) \tms_{\mathscr{R}^{\tau}(\Z)} \mathscr{R}^{\tau}_{irr}(K_2) \rt)     \sqcup \lt( \mathscr{R}^{\tau}_{irr}(K_1) \tms_{\mathscr{R}^{\tau}(\Z)} \mathscr{R}^{\tau}_{irr}(K_2) \rt). \nonumber
\end{align}

\lem \label{lemma:decomp-types}
	Suppose that $\tau \in \mc{G}(K_1) \cap \mc{G}(K_2) \sub \C- \{ \pm 2\}$. Then the open subscheme of irreducibles $\mathscr{R}^{\tau}_{irr}(K) \sub \mathscr{R}^{\tau}(K)$ consists precisely of the union of the second, third and fourth components in the above decomposition.  
\elem 

\pf 
It's clear that the second, third and fourth components consist of irreducible representations. Hence we only need to show that the first component does not contain an irreducible representations. Equivalently, we need to argue that an irreducible representation of $K$ cannot restrict to an abelian representation on both $K_1$ and $K_2$. This property was proved in \cite[Prop.\ 7.3]{cote-man}.
\epf

Following \cite[Def.\ 7.4]{cote-man}, we introduce the following definition: 

\defi 
	Let $K= K_1 \# K_2$ be as above. An irreducible representation $\rho: \pi_1(Y-K) \to \op{SL}(2,\C)$ is said to be of Type I if it restricts to an irreducible representation on $K_i$ and to an abelian representation on $K_j$ for $i, j \in \{ 1,2\}$, $i \neq j$.  An irreducible representation is said to be of Type II if it restricts to an irreducible representation on both factors. 
	We also refer to a connected component of $\mathscr{R}^{\tau}_{irr}(K)$ or $\mathscr{X}^{\tau}_{irr}(K)$ as being of Type I or Type II if its closed points are all of Type I or Type II respectively. So for instance, in the decomposition \eqref{equation:decomp}, the second and third terms are of Type I and the fourth term is of Type II.  
\edefi

The Type I locus of $\mathscr{X}^{\tau}(K)$ admits the following description:

\lem The image of the Type I locus $\lt( \mathscr{R}^{\tau}_{irr}(K_1) \tms_{\mathscr{R}^{\tau}(\Z)} \mathscr{R}^{\tau}_{ab}(K_2) \rt) \sqcup  \lt( \mathscr{R}^{\tau}_{ab}(K_1) \tms_{\mathscr{R}^{\tau}(\Z)} \mathscr{R}^{\tau}_{irr}(K_2) \rt)$ under the projection map $\pi:  \mathscr{R}^{\tau}(K) \to \mathscr{X}^{\tau}(K)$ is isomorphic to the disjoint union of $\mathscr{X}^{\tau}_{irr}(K_1)$ and $\mathscr{X}^{\tau}_{irr}(K_2)$. 
\elem

\pf 
	We only show that the image of $  \mathscr{R}^{\tau}_{ab}(K_1) \tms_{\mathscr{R}^{\tau}(\Z)} \mathscr{R}^{\tau}_{irr}(K_2)$ is isomorphic to $\mathscr{X}^{\tau}_{irr}(K_2)$ since the other case is analogous. According to \Cref{lemma:disjoint-union} and \eqref{equation:fiberconnect}, the map $\mathscr{R}^{\tau}(\Z) \to \mathscr{R}^{\tau}_{ab}(K_1)$ is induced by the composition $\pi_1(S^2-p-q) \to \pi_1(Y_1-K_1) \to H_1(Y_1-K_1;\Z) \simeq \Z$, which is an isomorphism.  The desired claim now reduces to a straightforward algebraic fact: let $A, B, B'$ be $\C$-algebras and suppose that $\op{SL}(2,\C)$ acts on the underlying vector spaces. Given a morphism $B \to A$ and an isomorphism $B \to B'$ which both commute with the $\op{SL}(2,\C)$ action, there is an induced action of $\op{SL}(2,\C)$ on the tensor product $A \otimes_B B'$ and an isomorphism of invariant rings $(A \otimes_B B')^{\op{SL}(2,\C)} \to A^{\op{SL}(2,\C)}$.  \epf

\subsection{A holomorphic $\C^*$-action} 

Let $(\S, p_i, q_i, U_i, U_i')$ be a Heegaard splitting for $Y_i-K_i$.  Following \cite[Sec.\ 7.4]{cote-man}, we can construct a Heegaard splitting for $Y-K= (Y_1 \# Y_2) - (K_1 \# K_2)$ as follows. Let $D_q \sub \S$ and $D_{q_1} \sub \S'$ be open balls around $q_1, p_2$ having smooth boundary and closure diffeomorphic to the unit disk. Fix a diffeomorphism $\psi: D_{q_1} \to D_{p_2}$ which extends smoothly to the boundary. We let $U_1 \#_b U_2$ and $U_1' \#_b U_2'$ be obtained by gluing the handlebodies via $\psi$. Let $\S \# \S'$ be obtained by gluing $\S-D_{q_1}$ and $\S-D_{p_2}$.  Observe that $\d(D_{q_1}) =_{\psi} \d(D_{p_2}) \sub \S \# \S'$ is a separating, simple closed curve which we call $c$. 

The goal of this section is to establish the following proposition. 

\prop \label{proposition:free-action}
For $\tau \in \mc{G}(K_1) \cap \mc{G}(K_2) \sub \C - \{\pm 2\}$, there exists a holomorphic action of $\C^*$ on an open subset $\mathcal{U} \subset \mathscr{X}_{irr}^{\tau} (\Sigma)$ such that the Lagrangians $L_0$ and $L_1$ are contained in $\mathcal{U}$ and preserved by the action. The induced action on the Type II locus of $\mathscr{X}_{irr}^{\tau} (K) =L_0 \cap L_1$ is free.
\eprop

The action which we will exhibit was already considered in \cite[Sec.\ 7.4]{cote-man}, but it will be useful to give a more detailed construction following \cite{goldman}. 

\subsubsection{A holomorphic action of $(\C, +)$} 

We assume throughout this section that $\tau \in \mc{G}(K_1) \cap \mc{G}(K_2) \sub \C - \{ \pm 2\}$. We begin with the following lemma. 

\lem
There exists a holomorphic, $\op{SL}(2,\C)$-equivariant action of the additive group $(\C, +)$ on $\mathscr{R}_{irr}^{\tau} (\Sigma)$, which therefore induces a holomorphic action on $\mathscr{X}^{\tau}_{irr}(\S)$.
\elem

\pf From van Kampen's theorem, we have the following description of $\pi_1(\Sigma)$:
\begin{align*}
\pi_1(\Sigma-\{p_1,q_2\}) =   \pi_1(\Sigma_1-\{p_1,q_1\})*_{\pi_1(\d D_{q_1})} \pi_1(\Sigma_2-\{p_2,q_2\}).
\end{align*}
We fix an isomorphism $\Z = \pi_1(\d D_{q_1})$ be sending $1 \mapsto c$. For $i=1,2$, we let $\iota_i: \Z \to (\Sigma_i-\{p_i,q_i\})$ be the maps inducing above pushout diagram. 

Identifying $c$ with its image under $\iota_i$, note that $c$ is a meridian for $K_i$.  The points of $\mathscr{R}^{\tau}(\S)$ can therefore be viewed as pairs $(\rho_1, \rho_2) \in \mathscr{R}^{\tau}(\S_1) \tms \mathscr{R}^{\tau}(\S_2)$ such that $\rho_1(c)=\rho_2(c)$.

Let $F: \op{SL}(2,\C) \to \mathfrak{sl}(2,\C)$ be the projection onto the trace-free part. That is, $F(A)=A-\frac{1}{2}\op{tr}(A)I$ for $A \in \op{SL}(2,\C)$. The additive group $(\C, +)$ acts on $\mathscr{R}^{\tau}(\S)$ by
\begin{align*}
t * (\rho_1,\rho_2)=\left( \exp(tF(\rho_1([\mu_1])))\rho_1\exp(-tF(\rho_1([\mu_1]))),\rho_2\right).
\end{align*}

This action is evidently holomorphic. It is well-defined on $\mathscr{R}^{\tau}(\Sigma)$ due to the fact that $\exp(t F(\rho_1(c))\in\op{Stab}(\rho_1(c))$ for all $t \in \C$. 

We claim that the action also restricts to $\mathscr{R}_{irr}^{\tau} (\Sigma)$. To prove this, it suffices to check that it sends reducibles to reducibles. If $(\rho_1,\rho_2)$ is reducible, then let $0\neq v\in\C^2$ be a generator of the line preserved by this representation, i.e. $v$ is an eigenvector of every matrix in the image. In particular, it is an eigenvector of $\rho_1(c)$, which means it is an eigenvector of $\exp(tF(\rho_1(c)))$. Thus, the line is also preserved by the representation $\exp(tF(\rho_1(c)))\rho_1\exp(-tF(\rho_1(c)))$.

The $\op{SL}(2,\C)$-equivariance of the action follows from the conjugation equivariance of $\exp$ and the projection $F$. It follows from the equivariance of the action that it passes to the quotient $\mathscr{X}_{irr}^{\tau} (\Sigma)$. \epf

\subsubsection{The induced $\C^*$ action}  Given a representation $\rho \in \mathscr{R}^{\tau}_{irr}(\S)$, the function $\zeta(\rho):=\op{det} F( \rho(c))$ is clearly algebraic and invariant under conjugation. It therefore defines a function on $\mathscr{X}^{\tau}_{irr}(\S)$ (i.e. an element of the ring of functions of this scheme). 
 
It's straightforward to check that $\zeta(\rho) \in \C^*$ due to our assumption that $\tau \neq \pm 2$. We can therefore choose a ball $B_{\e} \sub \C^*$ centered at $\zeta(\rho) \in \C^*$ and let $\mc{U}:= \zeta^{-1}(B_{\e}) \sub \mathscr{X}^{\tau}_{irr}(\S)$.  We also choose a square root on $B_{\e} \sub \C^*$ which will be fixed for the remainder of this section. 

Observe that the $(\C, +)$ action described in the previous section preserves $\mc{U} \sub \mathscr{X}^{\tau}_{irr}(\S)$. Let us now analyze the stabilizer of its restriction to $\mc{U}$. 

\lem \label{lemma:stabilizers}
	Given $[\rho] \in \mc{U}$, the stabilizer of $[\rho]$ under the $(\C, +)$ action contains $(\pi / \sqrt{ \zeta( \rho)}) \Z$. If $\rho = (\rho_1, \rho_2) \in \mathscr{R}^{\tau}(\S_1) \tms \mathscr{R}^{\tau}(\S_2)$ and $\rho_1, \rho_2$ are both irreducible, then the stabilizer is exactly $(\pi / \sqrt{ \zeta( \rho)}) \Z$.
\elem

\pf 
	By direct computation, one checks that $\op{exp}(t F( \rho(c))) = \pm \op{Id}$ if and only if $t \in (\pi / \sqrt{ \zeta( \rho)}) \Z$. 

	As we noted at the beginning of \Cref{subsection:structure-character-variety}, the irreducibility of $\rho_1$ and $\rho_2$ implies that their stabilizer under the conjugation action of $\op{SL}(2,\C)$ is $\pm \op{Id}$. Suppose now that $t*[(\rho_1, \rho_2)] = [t^*(\rho_1, \rho_2)] = [(\rho_1, \rho_2)]$. It follows by irreducibility of $\rho_2$ that $t* (\rho_1, \rho_2) = (\rho_1, \rho_2)$.  This then implies, by irreducibility of $\rho_1$, that $t \in (\pi / \sqrt{ \zeta( \rho)}) \Z$. \epf

\cor
	The restriction of the $(\C, +)$ action to $\mc{U}$ induces a holomorphic $\C^*$ action. 
\ecor

\pf
	
Given $\lambda \in \C^*$, define
\begin{align*}
\lambda \cdot [\rho] = \frac{1}{2 \sqrt{\zeta(\rho)}} \log(\lambda) * [\rho].
\end{align*}
It follows from the previous lemma that this action is well-defined (i.e. independent of the choice of logarithm). \epf

\pf[Proof of \Cref{proposition:free-action}] 
	We will treat the case of $L_0$ as the other one is analogous. If $\rho \in L_i$, then $\op{Tr} \rho(c) = \tau$ because $c$ is a meridian of $K$. It follows that $\rho \in \mc{U}$.  To see that the action preserves $L_0$, observe that a representation $\rho \in L_0$ can be viewed as a pair $\rho= (\rho_1, \rho_2)$ where $\rho_i \in \mathscr{R}^{\tau}(U_i)$ and $\rho_1(c)= \rho_2(c)$. Evidently, for $t * (\rho_1, \rho_2)$ is of the same form for $t \in \C$ so the claim follows.  Finally, the fact that the action is free on the Type II locus of $\mathscr{X}^{\tau}_{irr}(K)$ is an immediate consequence of \Cref{lemma:stabilizers}. \epf

Let us now consider the inclusion $\Z/n \hookrightarrow \C^* $ sending $[k] \mapsto e^{2\pi i k/n}$. It follows from \Cref{proposition:free-action} that the $\C^*$ action we have described is free on the Type II locus.  It follows that the induced $\Z/n$ action is free.

\cor \label{corollary:typeiizero}  The Euler characteristic of the Type II locus, weighted by the Behrend function, is zero. \ecor

\pf

Let us write $X_{\op{II}} \sub \mathscr{X}^{\tau}_{irr}(K_1 \# K_2)$ for the Type II locus. By definition we have $\chi_B(X_{\op{II}})= \sum_{m \in \Z} m \chi( \nu^{-1}(m)  )$. According to \cite[Prop.\ 4.2]{joyce-song}, the Behrend function depends only on the complex-analytic structure and is in particular preserved under isomorphisms of complex-analytic spaces. Hence, the $\Z/n$-action preserves the sets $X_{\op{II}}(m):= \nu^{-1}(m) $; in particular, $\Z/n$ acts on each $X_{\op{II}}(m)$.

The $\Z/n$ action on $X_{\op{II}}(m)$ is evidently free and properly discontinuous. Hence the quotient projection $X_{\op{II}}(m) \to X_{\op{II}}(m)/(\Z/n)$ is a covering map. 
Note that $X_{\op{II}}(m)= \nu^{-1}(m)$ is a pre-stratified subset of $\C^n$, in the sense of \cite{matherold}.  It follows by the main result of \cite{goresky} that $X_{\op{II}}(m)$ admits a triangulation of dimension at most $n$. Hence $X_{\op{II}}(m)$ is naturally a countable, locally finite CW complex. It then follows by an argument due to Belegradek \cite{belegradek} that the quotient $X_{\op{II}}(m)/(\Z/n)$ is homotopy-equivalent to a CW complex. Hence by \cite[Sec.\ 5.1]{mccleary}, we have a Serre spectral sequence for homology associated to the fibration $(\Z/n) \to X_{\op{II}}(m) \to X_{\op{II}}(m)/ (\Z/n)$. It follows from the existence of this spectral sequence that $\chi(X_{\op{II}}(m))$ is divisible by $\chi(\Z/n)=n$ for every $n$. Hence $\chi(X_{\op{II}}(m))=0$. Hence $\chi_B(X_{\op{II}})=0$. 
\epf
We now have the necessary tools to prove \Cref{theorem:connectedsum}.
\pf[Proof of \Cref{theorem:connectedsum}] We can assume that $\tau \in \mc{G}(K_1) \cap \mc{G}(K_2) \sub \C - \{\pm 2\}$. Note that the quotient map $\mathscr{R}^{\tau}(K) \to \mathscr{X}^{\tau}(K)$ preserves the decomposition \eqref{equation:decomp}.  It now follows from \Cref{lemma:decomp-types} that we have a scheme theoretic decomposition $$\mathscr{X}_{irr}^{\tau}(K)= \mathscr{X}^{\tau}_{irr}(K_1) \sqcup \mathscr{X}^{\tau}_{irr}(K_2) \sqcup X_{\op{II}}.$$

It follows from the definition of $\chi_B(-)$ that it is additive under disjoint unions of schemes. Hence 
\eqs \chi_B(\mathscr{X}^{\tau}_{irr}(K)) = \chi_B( \mathscr{X}^{\tau}_{irr}(K_1)) +\chi_B(  \mathscr{X}^{\tau}_{irr}(K_2) ) + \chi_B(X_{\op{II}})  = \chi_B( \mathscr{X}^{\tau}_{irr}(K_1)) +\chi_B(  \mathscr{X}^{\tau}_{irr}(K_2), \eeqs where we have used \Cref{corollary:typeiizero}.  The desired conclusion now follows from \eqref{equation:behrend-equality}. 
\epf

\section{The $\op{SL}(2,\C)$-Floer homology of a knot in families} \label{section:complexplane}

In this section, we prove that the $\tau$-weighted sheaf-theoretic Floer homology groups $\HP^*_{\tau}(-)$ constructed for $\tau \in (-2,2)$ by Manolescu and the first author in \cite{cote-man} can in fact be defined for all $\tau \in \C - \{ \pm 2\}$. As discussed in the introduction, the statement of \Cref{theorem:cassonconstant} implicitly relies on this fact. As another application, we observe that the groups $\HP^*_{\tau}(K)$ are canonically the stalks of a constructible sheaf $\mc{F}(K) \in D^b(\C- \{ \pm 2\})$ for a wide class of knots $K \sub Y$. We describe $\mc{F}(K)$ explicitly when $K \sub S^3$ is the figure-eight knot.

\subsection{Overview of the argument} \label{subsection:overview-complex-plane} Given data $K \sub Y$ and a choice of Heegaard splitting $(\S, U_0, U_1)$ for the knot exterior $E_K$ as in \Cref{subsection:construction}, the construction of $P^{\bullet}_{L_0, L_1}$ in \cite[Sec.\ 3]{cote-man} works for all $\tau \in \C - \{ \pm 2\}$. (When $\tau = \pm 2$, the arguments used to show that the Heegaard splitting gives smooth symplectic manifolds and smooth Lagrangians break down; see for instance Propositions 3.3 and 3.4 in \cite{cote-man}.)  The only place where one uses the assumption that $\tau \in (-2,2)$ is in proving that $P^{\bullet}_{L_0, L_1}$ is independent of choice of Heegaard splitting. This is done in \cite[Prop.\ 3.9]{cote-man}, where one crucially needs the fact that $X^{\tau}_{irr}(\S)$ is connected and simply-connected if $\S$ has genus at least $6$.  This fact is established in the appendix of \cite{cote-man} by exploiting a correspondence between character varieties of punctured surfaces and appropriate moduli spaces of parabolic Higgs bundles, whose topology is easier to analyze.  

Most of this section is devoted to proving that $X^{\tau}_{irr}(\S)$ is connected and simply-connected for all $\tau \in \C$ under the same assumption that $\S$ has genus at least $6$. This is the content of \Cref{proposition:candsc}. As explained above, it then follows immediately from the arguments in \cite{cote-man} that the perverse sheaf $P^{\bullet}_{\tau}(K)$ is well-defined for all $\tau \in \C- \{ \pm 2\}$; see \Cref{corollary:welldef}.  

The proof begins with the observation that the varieties $X^{\tau}_{irr}(\S)$ are the fibers of the projection map $X_{irr}(\S) \to \C$ described in \Cref{subsection:construction} which takes a representation to the trace of its holonomy along the boundary circles. It was already proved in \cite{cote-man} that $X^{\tau}_{irr}(\S)$ is connected and simply-connected for $\tau \in (-2,2)$. The key step is then to appeal to a theorem of Verdier (\Cref{theorem:verdier}) which implies that all but finitely many fibers are homeomorphic. It particular, all but finitely many fibers are connected and simply-connected.  

The conclusion can be extended to the remaining finitely many fibers by essentially repeating the arguments of the appendix of \cite{cote-man} in a slightly more general setting. Specifically, one exploits the correspondence between $X^{\tau}_{irr}(\S)$ and the moduli space of so-called $K(D)$-pairs, which are a slight generalization of the parabolic Higgs bundles considered in \cite{cote-man}.  These moduli spaces depend on a choice of ``weight data", which determine the stability conditions but do not affect the underlying bundles. One then shows by varying the weights that each fiber is homeomorphic to infinitely many other ones, and the conclusion follows. (For a few special cases, one needs to be more careful, as one only gets a homeomorphism in the complement of a set of large codimension. This sort of phenomenon also occurred in \cite{cote-man}.)

\subsection{Application of a Theorem of Verdier}
Let $\mc{S}$ be a compact and Riemann surface with empty boundary of genus $g \geq 2$ and let $p, q \in \mc{S}$ be a pair of distinct points.  Let $\S:= \mc{S}- p-q$. Let $c_p$ and $c_q$ be loops around $p$ and $q$ respectively, with respect to some arbitrary fixed basepoint. 

We choose a presentation $$\pi_1(\S)= \{ a_1,\dots,a_g,b_1,\dots,b_g,c_1,c_2 \mid \prod_{i}[a_i,b_i]c_1c_2=\op{Id} \}.$$ 

Let $\mathscr{R}(\S)$ be the representation scheme of $\mc{S}$. Using the above presentation for $\pi_1(\S)$, the ring of functions of $\mathscr{R}(\S)$ can be constructed as follows (see \cite[Sec.\ 2.1]{cote-man}): start with the polynomial ring in $2g+2$ variables $$\C[x^{a_1}_{11}, x^{a_1}_{12}, x^{a_1}_{21}, x^{a_1}_{22},\dots,x^{c_2}_{11}, x^{c_2}_{12}, x^{c_2}_{21}, x^{c_2}_{22}],$$  then mod out by the relations $\{x^{\a}_{11}x^{\a}_{21}-x^{\a}_{12}x^{\a}_{22}=1\}$ where $\a$ ranges over the generators of $\pi_1(\S)$. Finally, mod out by an additional polynomial coming from the relation $\prod_{i}[a_i,b_i]c_1c_2=\op{Id}$. 

Let $A$ be the resulting ring. The complex algebraic group $\op{SL}_2$ acts by conjugation. Let us call the ring of invariants $A^G$. 

Using the fact that $\op{SL}_2$ is linearly reductive, the coordinate ring of $\mathscr{X}^{\tau}(\S)$ is exactly 
\begin{align*} A^G/(x^{c_1}_{11}+ x^{c_1}_{22} - \tau, x^{c_2}_{11}+ x^{c_2}_{22} - \tau) &= A^G/(x^{c_1}_{11}+ x^{c_1}_{22} - \tau, x^{c_1}_{11}+ x^{c_1}_{22} - (x^{c_2}_{11}+ x^{c_2}_{22})) \\
&= A^G/ (x^{c_1}_{11}+ x^{c_1}_{22} - (x^{c_2}_{11}+ x^{c_2}_{22})) \otimes_{\C[x]} \C[x]/(x-\tau),
\end{align*} where the tensor product is formed via the map $\C[x] \to A$ sending $x \mapsto x^{c_1}_{11}+ x^{c_1}_{22}$. 

On the other hand, using again the fact that $\op{SL}_2$ is linearly reductive, the ring $A^G/ (x^{c_1}_{11}+ x^{c_1}_{22} - (x^{c_2}_{11}+ x^{c_2}_{22}))$ is the ring of invariants of $A/ (x^{c_1}_{11}+ x^{c_1}_{22} - (x^{c_2}_{11}+ x^{c_2}_{22}))$. Let $$W= \op{Spec} \lt(A^G/ (x^{c_1}_{11}+ x^{c_1}_{22} - (x^{c_2}_{11}+ x^{c_2}_{22})) \rt).$$  There is an open locus $W_{irr} \sub W$ corresponding to irreducible representations. 

\lem We have $\mathscr{X}^{\tau}_{irr}(\mc{S})= W_{irr} \times_{\bb{A}^1} \{\tau\} := (W_{irr})_{\tau}.$ \elem
\pf As demonstrated above, the coordinate rings of $\mathscr{X}^{\tau}(\mc{S})$ and $W\times_{\bb{A}^1}\{\tau\}$ are identical, so they are the same scheme. The irreducible representations in either scheme are precisely those that are irreducible as representations of $\pi_1(\mc{S})$. That is, $W_{irr}=W\cap \mathscr{X}_{irr}(\mc{S})$. Taking fibers over $\tau$, we find that $(W_{irr})_{\tau}=W_{\tau}\cap \mathscr{X}_{irr}(\mc{S})=\mathscr{X}^{\tau}(\mc{S})\cap \mathscr{X}_{irr}(\mc{S}) = \mathscr{X}_{irr}^{\tau}(\mc{S})$. \epf

We will need the following result of Verdier.

\thm[\cite{verdier}, Cor.\ 5.1] \label{theorem:verdier} Let $X, Y$ be complex algebraic varieties (separated and of finite type). Let $f: X \to Y$ be a morphism. Then there is a dense Zariski open set $U \sub Y$ such that $f|_U:= f^{-1}(U) \to U$ is a locally trivial topological fibration (in the analytic topology). \ethm 

We are led to the following corollary. 

\cor \label{corollary:scdense} There is a Zariski open set $\mc{V} \sub \bb{A}^1$ such that $X^{\tau}_{irr}(\mc{S})$ is connected simply-connected, and of dimension $6g-2$ for all $\tau \in \mc{V}$. \ecor

\pf Since $W_{irr}$ is an open subset of a finitely-generated complex algebraic variety, it is separated and of finite type.  So the theorem implies that it is a locally trivial topological fibration over some Zariski open subset $\mc{V} \sub \C$. 

It was shown in \cite[Appendix I]{cote-man} that $X^{\tau}_{irr}(\mc{S})$ is connected, simply-connected and of dimension $6g-2$ for all $\tau \in (-2,2)$. Since $(-2,2) \cap \mc{V}$ must be nonempty, it follows by \Cref{theorem:verdier} that $X^{\tau}_{irr}(\mc{S})$ is connected and simply-connected for all $\tau \in U \sub \mc{C}$.  \epf

\subsection{Moduli spaces of $K(D)$-pairs} 

We will now upgrade \Cref{corollary:scdense} by proving that $X^{\tau}_{irr}(\S)$ is in fact connected and simply-connected for all $\tau \in \C - \{\pm{2}\}$.  The argument makes use of a diffeomorphism between $X^{\tau}_{irr}(\S)$ and certain moduli spaces of so-called $K(D)$-pairs, which are a modest generalization of Higgs bundles.  

We note that \Cref{corollary:scdense} already allows us to make sense of $\HP^*_{\tau}(-)$ for generic $\tau$, which is all that one needs for the purpose of \Cref{theorem:cassonconstant} and \Cref{theorem:connectedsum}. However, we have chosen to include this section for completeness and in view of the possibility of studying $\HP^*_{\tau}(-)$ in families that was alluded to in the introduction. 

Throughout this section, we need to appeal to the general theory of parabolic vector bundles and Higgs bundles. The relevant definitions are introduced in Section 8.1 of \cite{cote-man} and we have chosen not to repeat them here for the sake of concision.  We now introduce a class of objects which are very similar to Higgs bundles and were not considered in \cite{cote-man}. A good reference for these is \cite{mondello} (but the reader should be warned that the objects which we refer to as $K(D)$-pairs are just called ``parabolic Higgs bundles" in \cite{mondello}). 

\defi[cf.\ Sec.\ 8.2 of \cite{cote-man}]  Fix a Riemann surface $\mc{S}$ and let $D= p_1+\dots+p_n$ for some collection of $n$ distinct points. A $K(D)$-pair on $(\mc{S}, D)$ is the data of a pair $(E_*, \Phi)$ consisting of a parabolic vector bundle $E_*$ and a (not necessarily strongly) parabolic morphism $\Phi: E_* \to E_* \otimes K(D)$, where $K$ is the canonical bundle. The morphism $\Phi$ is often called a \emph{Higgs field}. We usually denote $K(D)$-pairs by boldface letters $\bf{E}=(E_*,\Phi)$. \edefi

We remind the reader that parabolic Higgs bundles are defined in the same way as $K(D)$-pairs, except that one requires $\Phi$ to be a strongly parabolic morphism; see \cite[Sec.\ 8.2]{cote-man}. In fact many authors including \cite{mondello} refer to $K(D)$-pairs as parabolic Higgs bundles.  There are many other inconsistent conventions in this theory, so we remind the reader that we will always follow the conventions of \cite[Sec.\ 8]{cote-man}. 

For the remainder of this section, we specialize to the case of a Riemann surface $\mc{S}$ and a divisor $D=p+q$ for two distinct points $p,q \in \mc{S}$.  Let $\bm{\o}$ denote the data of weights $0 \leq \a_1(p) \leq \a_2(p) < 1$ and $0\leq \a_1(q) \leq \a_2(q)< 1$. Let $\fk{c}$ denote the data of a pair of matrices $\nu_p, \nu_q \in \fk{sl}(2,\C)$. 

We consider the moduli space $\op{wHiggs}^s(\mc{S}, \bm{\o}, 2, \mc{O}_{\mc{S}}, \fk{c})$ which parametrizes isomorphism classes of stable $K(D)$-pairs $(E_*, \Phi)$ satisfying the following conditions:
\begin{itemize}
\item $(E_*, \Phi)$ has rank $2$, 
\item $\op{Res}_p \Phi = \nu_p$ and $\op{Res}_q \Phi = \nu_q$,
\item the weights are given by $\bm{\o}$,
\item $\op{det}(E_*) \simeq \mc{O}_\mc{S}$,
\item $\op{Tr} \Phi=0$.
\end{itemize}

We define $\op{wHiggs}^{ss}(\mc{S}, \bm{\o}, 2, \mc{O}_{\mc{S}}, \fk{c})$ analogously, though we warn the reader that this is not in general a fine moduli space. We remark that the notation $\op{wHiggs}(-)$ is intended to be compatible with the notation of \cite[Sec.\ 8.2]{cote-man} and \cite{mondello}. The prefix ``w" stands for ``weak" and reflects the fact that the Higgs field of a $K(D)$-pairs satisfies a weaker condition than for an ordinary Higgs bundle. 

For $\a \in (0,1/2)$, it will be convenient to let $\bm{\o}(\a)$ denote the data of weights $0< \a < 1-\a<1$ at $p$ and $q$. If $\a =0$, we let $\bm{\o}(\a)= \bm{\o}(0)$ denote the weights $0 = \a_1(p)= \a_2(p)= \a_1(q)=\a_2(q)$. If $\a=1/2$, we let $\bm{\o}(\a)= \bm{\o}(1/2)$ denote the weights $1/2 = \a_1(p)= \a_2(p)= \a_1(q)=\a_2(q)$.  For $t \in [0, \infty)$ we let $\fk{c}(t)$ be the data of weights $\nu_p= \nu_q$ having $t$ as an eigenvalue. 

For our purposes, the importance of $K(D)$-pairs is mainly due to the following theorem. 

\thm[see Thm.\ 4.12 in \cite{mondello} and c.f.\ Thm.\ 8.4. in \cite{cote-man}] \label{theorem:kdequivalence} For $\tau \in \C$, choose $0 \leq \a \leq 1/2$ and $t>0$ so that $\tau = \op{Tr} (\op{diag}( e^t e^{2\pi i \a}, e^{-t} e^{-2 \pi i \a} ) )= e^t e^{2\pi i \a}+ e^{-t} e^{-2 \pi i \a}$. Then there is a real-analytic diffeomorphism  $$X^{\tau}_{irr}(\S) \simeq \op{wHiggs}^{s}(\mc{S}, \bm{\o}(\a), 2, \mc{O}_{\mc{S}}, \fk{c}(t)).$$  
\ethm 

In the proof of the next proposition, it will be convenient to view the choice weights as an additional piece of data on a fixed $K(D)$-pair. From this perspective, when one changes the weights, one does not change the underlying set of $K(D)$-pairs but one changes their slopes; i.e. one changes the stability conditions. 

\prop \label{proposition:candsc}
	For all $\tau \in \C$, the relative character variety $X^{\tau}_{irr}(\S)$ is connected and simply-connected. 
\eprop

\pf[Proof of \Cref{proposition:candsc}] 

	According to \Cref{corollary:scdense}, the proposition is already proved for all values of $\tau$ contained in a (Zariski) open set $\mc{V} \sub \bb{A}^1$ which contains the interval $(-2,2)$. 
	
	Fix $\tau \in \C$ and choose $(t, \a) \in \R_{\geq 0} \tms [0,1/2]$ so that $e^t e^{2\pi i \a} + e^{-t} e^{- 2\pi i \a} = \tau$. Observe that there exists $ \a' \in (0,1/2)$ so that $\tau' = e^t e^{2\pi i \a'} + e^{-t} e^{- 2\pi i \a'} \in \mc{V}$. We now consider three possibilities.  
	
	Case I: $\a \in (0,1/2)$. Given a $K(D)$ pair $(E_*, \Phi)$ of rank $2$, parabolic degree $0$ and weights $\bm{\o}(\a)$, it's not hard to check that the stability conditions are constant under varying $\a \in (0,1/2)$. We can therefore define a map $$\psi_{\a}: \op{wHiggs}^{s}(\mc{S}, \bm{\o}(\a), 2, \mc{O}_{\mc{S}}, \fk{c}(t)) \to \op{wHiggs}^{s}(\mc{S}, \bm{\o}(\a'), 2, \mc{O}_{\mc{S}}, \fk{c}(t))$$ by sending $(E_*, \Phi)$ to itself and replacing the weights $(\a, \a)$ by $(\a', \a')$. Since this map is evidently invertible, it is an isomorphism. We conclude that the left hand side is connected and simply connected since the right hand side is.  
	
	Case II: $\a =0$. Fix $z \in \S- p -q$. We define a map $$\psi_0: \op{wHiggs}^{s}(\mc{S}, \bm{\o}(0), 2, \mc{O}_{\mc{S}}, \fk{c}(t)) \to \op{wHiggs}^{s}(\mc{S}, \bm{\o}(\a'), 2, \mc{O}_{\mc{S}}, \fk{c}(t))$$ by sending $(E_*, \Phi) \mapsto (E_* \otimes \mc{O}(-z), \Phi)$ and replacing the weights $(0, 0)$ by $(\a', \a')$. 
	
	This map is an embedding, but it fails to be surjective. Indeed, $(E_*, \Phi)$ may admit a sub-bundle $(E^-, \Phi^-)$ of parabolic degree $-1+2 \a<0$, having weights $0 <\a<1$ at $p,q$.  Such a bundle is not in the image of $\psi_0$, and one can easily check that these are the only bundles which can fail to be in the image of $\psi_0$.

	Let $B= \op{wHiggs}^s(\mc{S}, \bm{\o}(\a), 2, \mc{O}_{\mc{S}}, \fk{c}(t))- \op{im}(\psi_0)$. We just saw that $B$ is contained in the locus $\mc{E}$ of the moduli space which consists of extensions of $K(D)$ pairs $0 \to \bf{E}^- \to \bf{E} \to \bf{E}^+ \to 0$ where $E^-$ has parabolic degree $-1+ 2\a$. Since $\op{wHiggs}^s(\mc{S}, \bm{\o}(\a), 2, \mc{O}_{\mc{S}}, \fk{c}(t))$ has dimension $6g-2$ by \Cref{corollary:scdense}, it follows that  $\op{im}(\psi_0)$ also has this property whenever $\op{dim}(B) \leq 6g-5$.
		
	The dimension of the space of extensions $K(D)$ pairs can be computed as in \cite[Sec.\ 8.4]{cote-man}, so we only sketch the details. Let $\mc{E}$ be the space of extensions. Let $\mc{X}$ be the set of pairs $(E^-, \Phi^-), (E^+, \Phi^+)$ where $E^+, E^-$ have rank $1$ and the underlying line bundles have degree $-1$. There is a natural forgetful map $\mc{E} \to \mc{X}$ and the dimension of $\mc{E}$ is bounded above by the sum of the dimension of $\mc{X}$ and of the fibers. 
	
	The dimension of $\mc{X}$ is computed in \cite[p.\ 3]{boden-yok} to be $2g+1$. The fiber over a fixed pair $(E^-, \Phi^-), (E^+, \Phi^+)$ is the space of extensions of this pair. Since each pair has different weights, this space is parametrized by the first homology group of an appropriate double complex as in \cite[Prop.\ 8.12]{cote-man}. The dimension can be computed as in \cite[Sec.\ 8.4]{cote-man} and one gets an upper bound of $4g+1$ on the dimension of $\mc{E}$. (This is essentially the same answer as in \cite[Thm.\ 8.9]{cote-man}, up to an additive constant which is independent of $g$). In particular, our assumption that $g \geq 6$ implies that $\op{dim}(B) \leq \op{dim}(\mc{E}) \leq 6g-5$ as desired. (This is in fact true once $g \geq 4$, but the requirement that $g \geq 6$ is needed in \cite[Prop.\ 3.14]{cote-man}). 
		
	Case III: $\a=1/2$.  We use the same map as in Case I. The map is an embedding, but fails to be a surjective as in Case II. The problem occurs again with sub-bundle of parabolic degree $-1+ 2\a$, and the subsequent argument is then the same as in Case II. \epf

\cor  \label{corollary:welldef}
	Given a knot $K$ in an oriented, closed $3$-manifold $Y$, the perverse sheaf $P^{\bullet}_{\tau}(K)$ constructed in \cite{cote-man} is well-defined (i.e. independent of the choice of Heegaard splitting) for $\tau \in \C - \{ \pm 2\}$. 
\ecor 
\pf
As explained in \Cref{subsection:overview-complex-plane}, the proof is entirely similar to the construction in Section 3 of \cite{cote-man}.
\epf

\subsection{$\op{SL}(2,\C)$ Floer homology in families} \label{subsection:families}

As an consequence of \Cref{corollary:welldef}, it makes sense to study the behavior of $\HP^*_{\tau}(K)$ in families for a given knot $K \sub Y$.  As discussed in the introduction, one expects that these groups should restrict to a local system on a Zariski open subset of the complex plane.  This expectation can already be verified for a wide class of knots considered in \cite{cote-man}. 

More precisely, for a fixed $3$-manifold $Y$, one considers \cite[Sec.\ 5.2]{cote-man} the class of all knots $K \sub Y$ whose character scheme $\mathscr{X}(K)$ is reduced and of dimension at most 1; see Assumptions A.1 and A.2 in \cite[Sec.\ 5.2]{cote-man}. For $Y=S^3$, this includes all two-bridge knots, torus knots and many pretzel knots.   Letting now $\pi: \mathscr{X}(K) \to \bb{A}^1$ be the map taking a representation onto its trace, it follows from the discussion in Section 5.3 of \cite{cote-man} that there is an open set $U \sub \bb{A}^1$ such that $\pi$ restricts to a smooth and proper morphism on the preimage of $U$. (More precisely, one can take $U$ is the set of points $\tau \in \bb{A}^1$ satisfying Assumptions B.1--B.4, which is shown to be a cofinite set). 

Letting $X_{irr}(K)_{\tau}$ be the fiber of $\tau \in \bb{A}^1$ under $\pi$, there is then a canonical identification $\HP^*_{\tau}(K) = H^*(X_{irr}(K)_{\tau})$; see \cite[Cor.\ 5.8]{cote-man}.  It is well-known that the cohomology of a smooth a proper map forms a local system on the base, so we conclude that there is a local system $\mc{F}_1(K) \in D^b(U)$ with $F(K)_{\tau} = \HP^*_{\tau}(K)$.  

Let $E(K)= \C- \{\pm 2\} - U$ (a finite set) and define $\mc{F}_2(K) \in D^b(E(K))$ to be the unique sheaf whose stalk at each point $\tau \in E(K)$ is $\HP_{\tau}^*(K)$. We then find that $$\mc{F}(K):= j_! \mc{F}_1 \oplus i_* \mc{F}_2 \in D^b(\C-\{\pm 2\})$$ is a constructible sheaf whose stalks compute $\HP^*_{\tau}(K)$, where $i: E(K) \to \C$ and $j:U \to \C$ are the inclusion maps. 

As mentioned in \Cref{subsection:construction}, it would be interesting to construct $\mc{F}(K)$ systematically, and for all knots, as the pushforward of a suitable sheaf on the character scheme $\mathscr{X}_{irr}(K)$ under the projection map considered above. 

\subsubsection*{An example: the figure-eight}  
According to \cite{porti}, the character variety (which agrees with the character scheme) of the figure-eight knot is \eq X(4_1)= \{(x,y) \mid (y-2)(y^2- (x^2-1)y+ x^2-1) =0 \}, \eeq 
where $\{(y-2)=0\}$ is the component of reducible representations and $x$ is the trace of a meridian.  

Let $E(4_1)= \{ \pm 1, \pm \sqrt{5} \}$. For $\tau \in \C - E(4_1)$, we have 
\eq \mathscr{X}^{\tau}_{irr}(4_1)=X^{\tau}_{irr}(4_1)= \{ (\tau, y) \mid y^2 - (\tau^2-1) y+ (\tau^2-1)=0 \} \sub \C. \eeq

For $\t \in [0, 2\pi]$ and $0< \e \ll 1$, let us consider the loop $\tau(\t)= 1+ \e e^{i \t}$. The points of $X^{\tau(\t)}_{irr}(4_1) \sub \C$ move around as 
$\t$ goes from $0$ to $2 \pi$ and can be computed by the quadratic formula. Noting that $\tau(\t)^2-1 = 2 \e e^{i \t} + \e^2 e^{2 i \t} \sim 2 \e e^{i \t}$, we compute that the relevant roots are approximately 
\eqs \frac{- 2 \e e^{i \t}  \pm \sqrt{ 4 \e^2 e^{i \t} - 8 \e e^{i \t}} }{2} \sim - \e e^{i \t} \pm \sqrt{2 \e} i e^{i \t /2}. \eeqs

This implies that the points of $X^{\tau(\t)}_{irr}(4_1) \sub \C$ get interchanged as $\t$ varies from $0$ to $2 \pi$. It follows from the above discussion that the local system $\mc{F}_1(4_1) \in D^b(\C - \{ \pm 1, \pm \sqrt{5}\} )$ has fibers $\Z^2$ concentrated in degree $0$. Fixing a reference fiber, we have seen that the monodromy around $+1$ is given by the matrix
$\begin{pmatrix}
0 & 1 \\
1 & 0 
\end{pmatrix}$. By a similar argument, one can check that the monodromy around $ \{ -1, \pm \sqrt{5} \}$ is given by the same matrix up to conjugacy.

\section{Appendix}

As explained in \Cref{subsection:constructible-theory}, there are many definitions in the literature of the Euler characteristic \emph{weighted} by a constructible function. The purpose of this appendix is to prove that these definitions are all equivalent. This is the content of \Cref{corollary:equiv-euler}.

\subsection{A new convention} It will be convenient in this appendix to consider ordinary Whitney stratifications (which we will call b-stratifications), rather than the more restrictive w-stratifications defined in \Cref{definition:wstrat} which we considered in the other sections of the paper. The reader may note that we only defined w-stratifications in the complex algebraic category (in particular, the strata were required to be complex algebraic). In contrast, we will define b-stratifications in the smooth category.  

\rmk 
	Vedier's condition w) also makes sense in the smooth category, so we could also have defined w-stratifications in the smooth category in \Cref{definition:wstrat}.  However, since this appendix is the only part of the paper in which we want to consider non-algebraic stratifications, it seemed more natural to include the requirement that the strata be algebraic as part of our definition. 
\ermk

For convenience, we only consider b-stratifications of subsets of $\R^n$, for $n \geq 1$. 

\defi[see Sec.\ 1 in \cite{verdier}] Let $(M', M)$ be a pair of locally closed, smooth submanifolds of $\R^n$ for some $n \geq 1$ such that $M \cap M'= \emptyset$ and $y \in \ov{M} \cap M'$. We say that the pair $(M',M)$ satisfies Whitney's condition b) at $y$ if the following holds: for any sequence $(x_n, y_n) \in M \tms  M'$ such that
\begin{itemize}
\item $x_n \to y$ and $y_n \to y$,
\item the sequence of lines $\R (x_n-y_n)$ has a limit $L$ in $\bb{P}(\R^n)$,
\item the sequence of tangent planes $T_{M, x_n}$ has a limit $T$ in $\op{Grass}(\R^n)$, 
\end{itemize}
then we have $L \sub T$. 
\edefi

We now state our notion of a b-stratification. We emphasize that this is not the most general definition (for instance, could allow locally-finite strata), but it is sufficient for our purposes. 

\defi[cf. (2.1) in \cite{verdier} and Sec.\ 1.2 in \cite{goresky-macpherson}] \label{definition:bstrat} A \emph{b-stratification} of a subset $X \sub \R^n$ is a partition $X= \sqcup_{i=1}^n X_i$, where the $X_i \sub X$ are locally closed, smooth submanifolds, which satisfies the following axioms:
\begin{itemize}
\item[(i)] $X_i \cap X_j= \emptyset$ if $i \neq j$. 
\item[(ii)] If $\ov{X}_i \cap X_j \neq \emptyset$, then $X_j \sub \ov{X}_i$. (One gets the same notion using the analytic or Zariski topology.)
\item[(iii)] If $X_i \sub \ov{X}_j$ and $i \neq j$, then the pair $(X_j, X_i)$ satisfies the condition b). 
\end{itemize}
Recall that the w-stratifications introduced in \Cref{definition:wstrat} are in particular b-stratifications. \edefi

\subsection{Equality of Euler characteristics}

Let $X$ be a closed subset of $\R^n$. Let $X= \sqcup_{i=1}^n S_i$ be a b-stratification which we call $\mc{S}$.  We say that an arbitrary subset $C \sub X$ is $\mc{S}$-constructible if $C= \cup_{ i \in \S} S_i$ for some subset $\S \sub \{1,2,\dots,n\}$.  Let $D=D(C)= \op{max}_{i \in \S} \op{dim}(S_i)$ and let $d=d(C)= \op{min}_{i \in \S} \op{dim}(S_i)$. We say that $D-d$ is the \emph{length} of $C$.

For $j= d, d+1,\dots,D$, let $C_j= \cup_{\op{dim}(S_k)=j} S_k$.  It follows from the second axiom of \Cref{definition:bstrat} that the subspace topology on $C_j$ coincides with the natural topology on $C_j$ as a disjoint union of smooth manifolds $S_k$. 

\defi[see p.\ 41 in \cite{goresky-macpherson}] Fix a point $x$ in a b-stratified set $X \sub \R^n$. Let $S$ denote some stratum and let $T$ be a smooth submanifold which is transverse to every stratum of $X$, which intersects $S$ at $x$ and nowhere else, and such that $\op{dim} S + \op{dim} T = n$. Let $B_{\d}(x):= \{ z \mid \| x-z \|  \leq \delta \}$ with the distance measured in the standard Euclidean metric.

For $0< \delta \ll 1$, let $N(x):= T \cap X \cap B_{\delta}(x)$ and let $\ell k(x) := T \cap X \cap \d B_{\delta}(x)$.  For $\delta$ small enough, the homeomorphism type of these spaces is independent of $T, \delta$, and of the choice of $x$. Moreover, they are canonically b-stratified as transverse intersections of b-stratified spaces. \edefi

\lem \label{lemma:localtrivial}There is a closed neighborhood $C_j \sub U_j$ with $U_j \sub C$ and a locally trivial projection map $\pi: U_j \to C_j$. The fiber $F_j$ over a point $x \in C_j$ is naturally a subspace of $N(x)$ and is $\mc{S}'$-constructible, where $\mc{S}'$ is the induced Whitney stratification on $N(x)$.  \elem

\pf According to \cite[p.\ 41]{goresky-macpherson}, there is a closed neighborhood $\tilde{U}_j \sub X$ of $C_j$  and a locally trivial projection map $\tilde{\pi}_j: \tilde{U}_j \to C_j$ whose fibers are homeomorphic to $N(x)$ for $x \in C_j$. Moreover, this fibration is locally homeomorphic to $\R^j \tms N(x)$ by a stratification-preserving homeomorphism.  The lemma now follows simply by letting $U_j= \tilde{U}_j \cap C$.  \epf

\cor The projection $U_j \to C_j$ is a (weak) homotopy equivalence (and hence an ordinary homotopy equivalence, since these are all CW complexes).  \ecor

\prop \label{proposition:euler-add}
	Suppose that one of the following two hypotheses holds:
\begin{itemize}
\item[(i)] All the strata have vanishing Euler characteristic,
\item[(ii)] all the strata are even-dimensional. 
\end{itemize}
Then for any $\mc{S}$-constructible subset $C= \cup_{i \in \S} S_i$, we have \eq \label{equation:eulerequal} \chi(C)= \sum_{i  \in \S} \chi(S_i)= \sum_{i \in \S} \chi_c(S_i)= \chi_{c}(C), \eeq where $\chi_c(C):= \sum_{i \in \Z} (-1)^i H^i_c(C)$. (In case (i), all these numbers are all zero!) \eprop

\pf  We work by induction on the length of $C$. If $C$ has length $1$, then $C$ is a union of smooth manifolds, which have vanishing Euler characteristic in case (i) and are even dimensional in case (ii). The desired result then follows from Poincar\'{e} duality. 

Suppose now that the result has been proved for all $\mc{S}$-constructible sets $C$ of length $n-1=D(C)-d(C)=D-d$. 

Let $C'= \cup_{j>d} C_j$. Then $C= C_d\cup C'$. By \Cref{lemma:localtrivial}, $C_d$ has a closed neighborhood $U_j$ (in general non-compact) such that $U_d-C_d$ is a locally trivial fibration over $C_d$. Since $\op{int}(U_d) \cup \op{int}(C')= C$, we have a Mayer-Vietoris long exact sequence for singular homology \eqs \dots \to H_k(U_d-C_d) \to H_k(C_d) \oplus H_k(C') \to H_k(C) \to \dots \eeqs

Suppose that (i) holds. Then $\chi(C_d)=0$. Hence $\chi(U_d-C_d)=0$ since $U_d-C_d$ is a locally trivial fibration. Since $\chi(C')=0$ by induction hypothesis, we conclude that $\chi(C)=0$. This proves the first equality of \eqref{equation:eulerequal} in case (i). 

Suppose now that (ii) holds.  We first claim that $\chi(U_d-C_d)=0$. Indeed, according to \Cref{lemma:localtrivial}, $U_d-C_d$ is a locally trivial fibration and according to \cite{sullivan}, the fiber $F_d$ satisfies (i).  We have already shown that this implies that $\chi(F_d)=0$. It follows that $\chi(U_d-C_d)=0$. The desired claim follows again from Mayer-Vietoris. This proves the first equality of \eqref{equation:eulerequal} in case (ii). 

The second equality of \eqref{equation:eulerequal} is a direct consequence of Poincar\'{e} duality. The third equality can be proved by the same argument as the first. One now needs to use a version of Mayer-Vietoris for compactly-supported cohomology (see for instance \cite[Sec.\ 3.4]{mathew}), as well as a multiplicativity property for Euler characteristic with compact support of locally trivial fibrations which can be deduced from the argument of \cite[Cor.\ 2.5.5]{dimca} using \cite[Cor.\ 2.3.24]{dimca}. (We note that there is a typo in the description of the $E_2$ page of the spectral sequence of \cite[Cor.\ 2.3.24]{dimca} which should be $E_2^{p,q}= H_c^p(Y, R^qf_! \mc{F}^{\bullet})$).  \epf

\cor \label{corollary:equiv-euler} Let $X$ be a complex-algebraic variety which admits an embedding into $\C^n$.  Let $f: X \to \Z$ be a constructible function.  Then all notions of the Euler characteristic of $X$ weighted by $f$ mentioned in \Cref{subsection:constructible-theory} coincide. More precisely, we have:  
\eq \chi(X,f):= \sum_{n \in \Z} n\cdot \chi(\{f^{-1}(n) \} ) = \sum_{S \in \mc{S}}  f|_{S} \cdot\chi(S)= \sum_{S \in \mc{S}}  f|_{S} \cdot \chi_c(S) = \sum_{n \in \Z} n \cdot\chi_c(\{f^{-1}(n) \} ) .\eeq  \ecor

\pf Since $f$ is constructible, it follows from \Cref{definition:constructible-function} and the comment following it that we can choose a w-stratification $X= \sqcup_{i=1}^n X_i$ (which is hence a b-stratification) by subvarieties which we call $\mc{S}$, and such that the sets $\{ f^{-1}(m)\}_{m \in \Z}$ are $\mc{S}$-constructible. The second and fourth equalities then follow from \Cref{proposition:euler-add}. The third equality is a consequence of Poincar\'{e} duality, since the strata are even-dimensional manifolds.
\epf

\end{document}